\documentclass[11pt]{amsart}
\usepackage[utf8]{inputenc}
\usepackage[margin=1.3in]{geometry}
\usepackage[titletoc,title]{appendix}
\usepackage[dvipsnames]{xcolor}
\usepackage{amssymb,amsfonts,amsthm,amscd,stmaryrd}
\usepackage{hyperref}
\usepackage{amsmath}
\usepackage{fancyhdr}
\usepackage{sidecap}
\usepackage{booktabs}
\usepackage{graphicx}
\usepackage{float}
\usepackage{mathtools}
\usepackage{indentfirst}
\usepackage{mathrsfs}
\usepackage{listings}
\usepackage{amsbsy}
\usepackage{amsfonts}
\usepackage{latexsym}
\usepackage{amsopn}
\usepackage{amsxtra}
\usepackage{euscript}
\usepackage{ dsfont }
\usepackage{ esint }
  \usepackage{mathrsfs}
 \usepackage{latexsym}
 \usepackage{bbm}
 \usepackage{hyperref}
 \usepackage{mathtools}
 \usepackage[font=small,format=hang,labelfont={sf,bf}]{caption}
 \usepackage{epsfig}
 \usepackage{subfig}
 \usepackage{url}
 \usepackage{varioref}
 \usepackage{bm}
\numberwithin{equation}{section}
\theoremstyle{plain}
\begingroup
\newtheorem{theorem}{Theorem}[section]
\newtheorem*{theorem*}{Theorem}
\newtheorem{lemma}[theorem]{Lemma}
\newtheorem{proposition}[theorem]{Proposition}

\endgroup
\theoremstyle{definition}
\theoremstyle{remark}
\begingroup
\newtheorem{definition}[theorem]{Definition}
\newtheorem{remark}[theorem]{Remark}

\endgroup
\newcommand{\bd}{\partial}

\newcommand{\R}{\mathbb{R}}

\newcommand{\N}{\mathbb{N}}
\newcommand{\Z}{\mathbb{Z}}
\newcommand{\T}{\mathbb{T}}

\newcommand{\Sp}{\mathbb{S}}
\newcommand{\ud}{\,\textnormal{d}}
\newcommand{\dist}{\text{dist}}
\newcommand{\virg}[1]{``#1''}
\let\div\undefined
\newcommand{\div}{\textnormal{div}}
\newcommand{\udH}{\ud\mathcal H^{N-1}}
\let\dist\undefined
\newcommand{\dist}{\textnormal{dist}}
\newcommand{\sd}{\textnormal{sd}}

\newcommand{\e}{\varepsilon}
\newcommand{\la}{\langle}
\newcommand{\ra}{\rangle}
\newcommand{\D}{\Delta}
\newcommand{\n}{\nabla}
\newcommand{\HHH}{\mathrm{H}} 

\newcommand{\beq}{\begin{equation}}
\newcommand{\eeq}{\end{equation}}
\newcommand{\pa}{\partial}

\usepackage[maxbibnames=99]{biblatex} 
\title[Stability geometric flows]{Stability of the surface diffusion flow and volume-preserving mean curvature flow in the flat torus}

\author[Daniele De Gennaro]
 {Daniele De Gennaro}
 \address[Daniele De Gennaro]{CEREMADE department, Université Paris-Dauphine, 1 Place du Maréchal de Lattre de Tassigny, 75775, Paris, France, \& Dipartimento di Scienze Matematiche, Fisiche e Informatiche, Università degli Studi di Parma, Parco Area delle Scienze 53/A, 43124, Parma, Italy}
 \email[Daniele De Gennaro]{degennaro@ceremade.dauphine.fr}

 \author[Antonia Diana]
 {Antonia Diana}
 \address[Antonia Diana]{Scuola Superiore Meridionale, Largo San Marcellino 10, 80138, Naples, Italy}
 \email[Antonia Diana]{antonia.diana@unina.it}
 
 \author[Andrea Kubin]
 {Andrea Kubin}
 \address[Andrea Kubin]{Zentrum Mathematik-M7, Technische Universit\"at M\"unchen, Boltzmannstrasse 3, 85747, Garching, Germany}
 \email[Andrea Kubin]{andrea.kubin@tum.de}
 
 \author[Anna Kubin]
 {Anna Kubin}
 \address[Anna Kubin]{Dipartimento di Scienze Matematiche ``G.~L.~Lagrange'', Politecnico di Torino, Corso Duca degli Abruzzi, 24, 10129 Torino, Italy \& Institute of Analysis and Scientific Computing, Technische Universit\"at Wien, Wiedner Haupstrasse 8-10, 1040 Vienna, Austria}
 \email[Anna Kubin]{anna.kubin@polito.it \& anna.kubin@asc.tuwien.ac.at}

\date{}
\begin{document}
\maketitle
\begin{abstract}
We prove that, in the flat torus and in any dimension, the volume-preserving mean curvature flow and the surface diffusion flow, starting $C^{1,1}-$close to a strictly stable critical set of the perimeter $E$, exist for all times and converge to a translate of $E$ exponentially fast  as time goes to infinity. 
\noindent 
\vskip5pt
\noindent
\textsc{Keywords}: Geometric evolutions; variational methods; volume constraint.  
\vskip5pt
\noindent
\textsc{AMS subject classifications:}  
53E10; 53E40;
49Q20; 35B40; 37E35.
\end{abstract}

\tableofcontents

\section*{Introduction}

In this paper we establish global in time existence and convergence towards equilibrium of two physically relevant volume-preserving geometric motions, namely the volume-preserving mean curvature flow and the surface diffusion flow. 

On the one hand, the first one is the volume-preserving counterpart of the well-known mean curvature flow, and it is defined as a smooth evolution of sets $E_t$ 
governed by the law
\begin{equation}\label{evol law MCF}
    V_t=-\HHH_{E_t}+\bar \HHH_{E_t}\quad \text{on }\bd E_t,
\end{equation}
where $V_t$ and $\HHH_{E_t}$ are the outer normal velocity
and the mean curvature of $\bd E_t$, respectively, while $\bar \HHH_{E_t}=\fint_{\bd E_t}\HHH_{E_t}$.
The mean curvature flow is a famous evolution model, with far-reaching geometric and physical applications, which has a rich history dating back to its use in material science. One notable application is in physical systems involving multiple phases, such as the motion of grain boundaries in materials science, as first discussed by Mullins \cite{Mullins}.

On the other hand, the surface diffusion flow is a smooth flow of sets $E_t$ evolving according to the law
\begin{equation}\label{evol law surf diff}
    V_t=\Delta_{E_t} \HHH_{E_t}\quad \text{on }\bd E_t,
\end{equation}
where $\Delta_{E_t}$ denotes the Laplace-Beltrami operator on $\bd E_t$.
Similar to the mean curvature flow, the surface diffusion flow has important applications in material science, especially in physical systems with multiple phases.
It has been proposed in the physical literature by Mullins \cite{Mul} to model surface dynamics for phase interfaces when the evolution is governed by mass diffusion in the interface.

The volume preserving mean curvature flow  can be seen as a simplified, second-order version of the surface diffusion flow as both flows share several common properties. Indeed, from the evolution laws \eqref{evol law MCF} and \eqref{evol law surf diff} it follows that the volume of the evolving sets is preserved along the two flows, as can be easily seen from the following computation
\[  \dfrac{\ud}{\ud t}|E_t|=\int_{\bd E_t}V_t\udH =0, \]
the perimeter is decreasing, since the evolution \eqref{evol law MCF} satisfies
\[ \dfrac{\ud}{\ud t} P(E_t)=\int_{\bd E_t} V_t \HHH_{E_t}\udH  =\int_{\bd E_t} \left( \HHH_{E_t} - \bar  \HHH_{E_t}  \right)^2\udH    \le 0,  \]
and an integration by parts shows for \eqref{evol law surf diff} that
\[  \dfrac{\ud}{\ud t} P(E_t)=\int_{\bd E_t} V_t \HHH_{E_t}\udH = -\int_{\bd E_t} |\n \HHH_{E_t}|^2\udH\le 0. \]
Moreover, these two evolutions can be regarded (at least formally) as gradient flows of the perimeter according to suitable metrics. In particular, the mean curvature flow can be considered as (a volume preserving modification of) the $L^2$-gradient flow of the perimeter, while the surface diffusion can be interpreted as its $H^{-1}$-gradient flow.

In both cases, singularities may appear in a finite time even for initial smooth sets (see \cite{MaySim}), therefore in general only short-time existence results are available, see for instance \cite{ES, Hui} for the mean curvature flow and \cite{EMS} for the surface diffusion flow (see also \cite{GarckGob1} for the case of triple junction clusters).  Because of the (formal) gradient flow structure of the two flows, it is reasonable to expect that if the initial set  is sufficiently close to a stable  point (or a local minimizer) $E$ of the perimeter, then the flow exists for all times and asymptotically converges to $E$.  We refer to this property as \textit{dynamical stability.}  We will properly define the notion of stability  in Definition \ref{def strictly stable}, however we can summarize it as follows: stable sets are sets whose boundary has constant mean curvature and positive definite second variation of the perimeter (i.e., they are \virg{stable} for the perimeter functional). In this paper, we will focus on the flat torus $\mathbb T^N$, which is particularly interesting due to the  variety of possible limit points of the flows, 
namely periodic constant mean curvature hypersurfaces. In the Euclidean space  only unions of balls have constant mean curvature, whereas the flat torus admits a much broader range of such surfaces. However, a full characterization of constant mean curvature hypersurfaces in $\T^N$ is not available in any dimension.  In dimension $N=2$, the only sets with constant mean curvature are discs and stripes (also called lamellae), while  for $N\ge 3$ there exist many nontrivial examples, as stripes, cylinders and triply periodic surfaces known as gyroids.

The aforementioned approach of studying the dynamical stability of stable sets has  been used in many instances in the literature. Concerning the surface diffusion, this method was employed in \cite{AFJM,FusJulMor2D,FusJulMor3D}, where the authors considered the surface diffusion (also with an extra elastic term) and the Mullins-Sekerka flows in the $2,3$-dimensional flat torus (see also the survey \cite{DDMsurvey}) and proved the dynamical stability of stable sets.  It should be noted that the flows considered in these works include nonlocal terms, but their results also apply to the evolution driven solely by the perimeter energy. In the Euclidan setting, other results for the surface diffusion deal with the stability of balls \cite{EMS, LeCSim20, Whe}, infinite cylinders \cite{LeCSim13}, two-dimensional triple junctions \cite{GarItoKoh}, as well double bubbles \cite{AbeAraGar,GarckGob2} (see also~\cite{LeCSim20} for similar results in different settings).\\ 
Regarding the volume preserving mean curvature flow, recent progresses have been made in proving the dynamical stability of strictly stable sets in the $3$-dimensional flat torus \cite{Nii}, while  older results mainly concern convex sets, balls, or the $2$-dimensional setting.
The dynamical stability of balls has been proven in the Euclidean setting  under various hypoteses on the dimension or on the initial set in \cite{ES, Gag, Hui,Li}. {We refer also to \cite{PruSim}, where global existence and convergence results for a large class of geometric evolution laws have been considered, relying on the concept of $L^p$-maximal regularity for quasilinear parabolic equations.} Another interesting approach, up to now limited to the  mean curvature flow, deals with the long-time behaviour of weak solutions of the flow, in particular the so-called flat flows \cite{ATW, LS}. Flat flows are measure-theoretic weak solutions to the mean curvature flow arising as the limit of a discrete-in-time approximation,  
based on the minimizing movement scheme,  as the time-step parameter tends to 0.  
The exponential convergence of flat flows to balls has been proved in \cite{JulMorPonSpa} in $\R^2$, while in \cite{BelCasChaNov} the authors deal with the anisotropic and crystalline mean curvature flow in $\R^N, N\ge 2,$ and for convex initial data, showing the asymptotic convergence to a Wulff shape. Concerning the time-discrete flows, in \cite{MorPonSpa} and \cite{DegKuKu}, the asymptotic convergence to balls in $\R^N, N\ge 2,$ is shown in the classical and fractional settings, respectively. In \cite{DeKu}, two of the authors prove  the dynamical stability for the discrete flow of strictly stable sets in the flat torus of any dimension.

\smallskip 
In the present paper we are able to prove in any dimensions the dynamical stability of strictly stable sets in the flat torus both for the surface diffusion flow and the volume preserving mean curvature flow.
{By  assuming the initial set to be close in the $C^{1,1}$-topology to a strictly stable set, we obtain global existence and asymptotic convergence of both the flows to (a translated of) the underlying stable set. This is quite surprising for the surface diffusion flow, which is a fourth-order flow. Our main result is the following.}

\begin{theorem}\label{teo asymptotic}
Let $E\subset \T^N$ be a strictly stable set and let $E_0 =E_{u_0} \subset \T^N$ be the normal deformation of $E$ induced by $u_0 \in C^{1,1}(\bd E)$ (see Definition \ref{normaldeformation}) with $|E_0|=|E|$. There exists $\delta=\delta(E)>0$ such that if $\|u_0\|_{C^{1,1}(\bd E)}\le \delta$, then 
\begin{itemize}
    \item[(i)] the volume-preserving mean curvature flow $E_t$ starting from $E_0$ (defined in \eqref{defmcf}) exists smooth for all times $t\ge 0$, and $E_t \to E +\tau$ as $t \to \infty$, for some $\tau \in \T^N$, in $C^k$ for every $k\in \N$ exponentially fast;
    \item[(ii)] the surface diffusion flow $E_t$ starting from $E_0$ (defined in \eqref{defsdf}) exists smooth for all times $t\ge 0$, and $E_t \to E +\tau$ as $t \to \infty$, for some $\tau \in \T^N$, in $C^k$ for every $k\in \N$ exponentially fast.
\end{itemize}
Where with exponentially fast we mean that the sets $E_t$ can be written as normal deformations of $E+\tau$ induced by functions $u(\cdot,t) \in C^{\infty}(\bd E+\tau)$ such that 
\[\|u(\cdot,t)\|_{C^k(\bd E+\tau)}\le C_k e^{-C_k t} \quad \text{for} \quad t>0.\]
\end{theorem}

{
The  main technical novelty of our argument is the use  a quantitative Alexandrov-type inequality, which has been obtained by two of the authors in \cite[Theorem~1.3]{DeKu} and is applied for the first time to a continuous-in-time setting, in this paper. This technique allows us to treat in a unified fashion both the geometric flows considered. 
However, it seems to be quite general, in the sense that it can be adapted to other gradient flows of the perimeter functional. 
For instance, we are confident that the  Mullins-Sekerka flow or, more in general, fractional gradient flows of the perimeter could be treated analogously,
provided one has sufficient control on the Schauder estimates for the linearized system governing the evolutions.  This will be the subject of future investigations. 
Moreover, since this stability inequality can be seen as a {\L}ojasiewicz-Simon inequality with sharp exponents, one is able to derive the optimal decay of the dissipation along the flow, immediately yielding the exponential convergence in any norm of the flow to the subjacent strictly stable set. 
In particular, our line of proof works in any dimension without the need of deriving energy estimates for the high derivatives of the curvature, which was one the main bottleneck of the previous methods developed in \cite{AFJM,FusJulMor2D,FusJulMor3D}.
Lastly, the Schauder-type estimates we provide following the lines of \cite{HZ} seems to be new in this setting.
}

We now outline the strategy of the proof, which is based on the gradient flow structure of the evolution.
Firstly, applying the  Alexandrov-type inequality  \cite[Theorem 1.3]{DeKu}, {combined with the quantitative isoperimetric inequality of \cite{AFM}}, we are able to to bound the velocity in terms of the displacement.
By iterating this procedure for the whole time of existence and using higher order estimates, we can extend the flow for all times.
In order to do so, we need to show that the short-time existence and regularity results depend only on the bounds of the initial datum. This is not \textit{a priori} clear from previous existence results \cite{EMS,ES}. More precisely, we rely on Schauder estimates on the linearized problem solved by the flows, which is a quasilinear perturbation of the heat equation for the mean curvature flow and a quasilinear perturbation of the biharmonic heat equation for the surface diffusion flow.
While Schauder-type estimates for general quasilinear parabolic PDEs of the second order are well known (see for instance \cite{Fri}), we couldn't find a precise reference for the fourth-order equation. Although an approach by scaling (in the spirit of \cite{KocLam}) could be feasible by working in local coordinates, we preferred to rely on the estimates provided in \cite{HZ}, where time-weighted H\"older norms are employed. 
After establishing the global existence of both flows, we obtain the exponential convergence up to translations via a Gronwall-type inequality.
This is where it comes into play the optimality of the exponent in the Alexandrov theorem \cite[Theorem 1.3]{DeKu}, which yields the exponential rate of convergence. 
Finally, we prove the convergence of these translations by exploiting the decay of geometric quantities along the flow, as in \cite{AFJM}.

{We conclude by highlighting that a similar stability result for the surface diffusion flow has been obtained by the second author and collaborators in \cite{DFM} using different techniques (that are shown in details in dimension $N=4$ and listed for any general $N$), assuming different hypotheses on the initial datum, depending on the dimension $N$. In particular, they consider initial sets $E_0$ close to the strictly stable set in $C^1$ and such that the  energy 
\[\int_{\bd E_0} |\nabla^{N-2}\mathrm H|^2  \ud \mathcal{H}^{N-1} + \int_{\bd E_0} |\nabla \mathrm H|^2  \ud \mathcal{H}^{N-1}\] 
is sufficiently small.
}

\section{Preliminary results}
In this section we collect some preliminary results and we fix the notations.

We denote by $\T^N$ the $N$-dimensional flat torus, which is the quotient of $\R^N$ by $\Z^N$.
The function spaces  $C^k(\T^N)$ and $W^{k,p}(\T^N)$, for $k \in \N$ and $p \in [1,\infty],$ are defined as the restriction of 
$C^k(\R^N)$  and $W_{loc}^{k,p}(\R^N)$, respectively, to the functions that are one-periodic. With $B_r(x)$ we denote the ball in $\R^N$ of center $x$ and radius $r$, while $B_r$ will be a short-hand notation for  $B_r(0)$. 
Given $x \in \R^N$, we will write $x=(x',x_N)$ where $x'\in \R^{N-1}$ and $x_N\in \R$.
Similarly, we denote by $B'_r(x')\subset \R^{N-1}$ the ball in $ \R^{N-1}$ with radius $r>0$ and center $x' \in \R^{N-1}$.

Moreover, we denote by $c,C $ some constants, which could be changing from line to line and always depend on the dimension $N$, and by $\frac {\bd }{\bd t}$ (or equivalently $\bd_t$) the partial derivative with respect to the variable $t$. 
{Let $F \subset \T^N$ we denote with $\dist_F(\cdot)$ the distance from the set $F$ and with $C^{1,1}(\pa F)$ the set of functions continuously differentiable with derivative Lipschitz continuous on $\pa F$.}

Given a smooth closed $(N-1)$-manifold $\Sigma \subset \T^N$ we denote by $\nu_{\Sigma}:\Sigma \to \Sp^N$ the outer normal to $\Sigma$, by $B_{\Sigma}$ the second fundamental form of $\Sigma$, and by $\HHH_{\Sigma}$ its mean curvature, that is the trace of $B_{\Sigma}.$ 
For every vector field $X:\Sigma \to \R^N$ we let $X_\tau$ to be the tangential part of $X$, that is $X_\tau(x)=X(x)-X(x) \cdot \nu_{\Sigma}(x) \nu_{\Sigma}(x)$, and for every function $f \in L^1(\Sigma)$ we denote with 
\[\bar f= \frac{1}{\mathcal{H}^{N-1}(\Sigma)}\int_{\Sigma} f d \mathcal{H}^{N-1}\]
the mean of $f$ over $\Sigma$.

Let $E \subset \T^N$ be a open set with smooth boundary and let $X:\T^N\to \R^N$ be a vector field of class $C^2$. We consider the associated flow $\Phi:\T^N\times (-1,1)\to\T^N$ defined by
\begin{equation}
  \frac {\bd \Phi}{\bd t}=X(\Phi) \quad \Phi(\cdot,0)=I, 
\end{equation}
where $I:\T^N \to \T^N$ denotes the identity,
and we say that $E_t =\Phi(E,t)$ is the variation of $E$ associated to $\Phi$ (or to $X$). If in addition it holds $\vert E_t \vert = \vert E \vert$ for every $t \in (-1,1)$, we say that $E_t$ is a \textit{volume-preserving variation} of $E$.

We now recall some  results on sets of finite perimeter, referring to \cite{Mag} for the basic definitions and proofs.  We say that a measurable set $E \subset \mathbb{T}^N$ is a set of finite perimeter if
 \begin{equation}
     P(E)\coloneqq \sup \left\{ \int_{E} \div(X)\ud x :  X \in C^1(\mathbb{T}^N,\R^N),\, \vert X \vert \leq 1 \right\} <\infty.
 \end{equation}
 Moreover, by De Giorgi's structure theorem, we have $P(E)= \mathcal{H}^{N-1}(\partial^* E)$ where $\partial^* E$ is a suitable $(N-1)$-rectifiable subset of $\partial E$. The {first and second variation of the perimeter at $E$} with respect to the flow $\Phi$ are defined as follows
\[\delta P(E)[X]\coloneqq \dfrac {\ud} {\ud t}\Big\lvert_{t=0}P(E_t), \quad \delta^2 P(E)[X]\coloneqq \dfrac {\ud^2} {\ud t^2}\Big\lvert_{t=0}P(E_t) .\]
It is well known  that, for any set of finite perimeter $E$, we have 
\begin{align*}
    \delta P(E)[X]=\int_{\partial^* E } \div_\tau(X)d\mathcal{H}^{N-1} ,
\end{align*}
where $\div_{\tau}(X)$ is the tangential divergence of $X$ on $E$ and it is given by
\begin{equation}
    \div_{\tau}(X)(x)= \div(X)(x)-\nu_{E}(x) \cdot \nabla X(x)  \quad \text{for all } x \in \bd^* E.
\end{equation}
Moreover, if $ E$ is a open set with $C^2$-boundary we have
\begin{align*}
    \delta P(E)[X]=\int_{\partial E } \div_\tau(X)d\mathcal{H}^{N-1} 
 =\int_{\bd E} \HHH_E \nu_E\cdot X\udH.
\end{align*}
Finally, the second variation formula for perimeter
on open sets of class $C^2$ (see for instance \cite[Section 3]{AFM}) is given by
\begin{align*}
		\delta^2 P(E)[X]=&\int_{\bd E} \left( |{\nabla_\tau} (X\cdot \nu_E)|^2-|B_E|^2 (X\cdot \nu_E)^2 \right)\udH \\
  &-\int_{\bd E} \HHH_E \div_\tau (X_\tau(X\cdot \nu_E)) \udH
		+\int_{\bd E} \HHH_E(\div X)(X\cdot \nu_E)  \udH,
\end{align*}
where {$\nabla_\tau f(x)= (\nabla f)_{\tau}(x)$} denotes the tangential derivative of $E$. 
Since the expression above only depends on the normal component of the velocity field $X$, we also denote by $\delta P(E)[\varphi]$ and $\delta^2 P(E)[\varphi]$, respectively, the first and the second variation of the perimeter at $E$, where $\varphi= X \cdot \nu_E$.

Let $E$ be a critical point of the perimeter. It should be noted that the translation invariance of the perimeter implies that the second variation becomes degenerate along flows of the form $\Phi(x,t)=x+t\eta,$ where 
$\eta\in \R^N$. Because of that, we denote by 
$$\tilde H^1(\bd E)\coloneqq \left\lbrace \varphi\in H^1(\bd E) : \int_{\bd E} \varphi \udH=0   \right\rbrace,$$
and by $T(\bd E)$ the subspace generated by the functions $\nu_i : \pa E \to \R$ for $i=1,\dots, N$, {defined as $\nu_i:=e_i \cdot \nu_E$ where $e_1, \dots, e_N$ is the standard orthonormal basis of $\R^N$.} We then set $T^\perp(\bd E)$ to be the orthogonal subspace of $T(\bd E)$ in the $L^2-$sense, that is
\[ T^\perp(\bd E)= \left\lbrace \varphi\in\tilde H^1(\bd E) : \int_{\bd E} \varphi  \nu_i\udH=0,\ i=1,\dots,N \right\rbrace.\]
After defining all the spaces, we can finally give the notion of stability.
\begin{definition}\label{def strictly stable}
We say that a set $E \subset \T^N$ of class $C^2$ is a \textit{strictly stable set} if it is a critical set, that is $\delta P(E)[\varphi]=0$ for all $\varphi \in \tilde H^1 (\bd E)$, and the second variation of the perimeter is positive definite, in the sense that 
	\[\delta^2P(E)[\varphi]>0,\qquad \forall \varphi \in T^{\perp} (\bd E) \setminus \{0\}.\]
\end{definition}

We now recall some technical results that will be useful in the following. 
We start by recalling the definition of normal deformation of a set and a result which ensures that any $W^{2,p}$-small normal deformation of a smooth set can be translated in a way so the projection on the subspace $T^\perp(E)$ becomes small.
\begin{definition}\label{normaldeformation}
Let $E \subset \mathbb{T}^{N}$ be an open set of class $C^{1}$. For every $f\in L^{\infty}(\bd E)$ such that $\|f \|_{L^{\infty}(\partial E)}$ is sufficiently small, we define the \textit{normal deformation of $E$ induced by $f$} the set $E_f$ having as boundary
$$\partial E_{f}: =\left\{x+ f(x)\nu_{E}(x) \, \colon \, x \in \partial E \right\}. $$
\end{definition}

\begin{lemma}[{\cite[Lemma 3.8]{AFM}}]\label{lemma 3.8_AFM}
Let $E\subset \T^N$ be of class $C^3$ and let $p>N-1$. For every $ \delta^*>0$ there exist $C>0$ and $\eta>0$ such that if $F$ is a normal deformation of $E$ induced by some $\psi\in C^2({\bd E})$ with $\|\psi\|_{W^{2,p}(\bd E)}\le \eta$, that is $F=E_\psi$, then there exist $\sigma\in\T^N$ and $\varphi\in W^{2,p}(\bd E)$ with the properties that 
\[ |\sigma|\le C\|\psi\|_{W^{2,p}(\bd E)},\quad \|\varphi\|_{W^{2,p}(\bd E)}\le C\|\psi\|_{W^{2,p}(\bd E)} \]
and 
\[ F+\sigma=E_\varphi, \quad \left\lvert \int_{\bd E}\varphi\nu_E\udH  \right\rvert\le  \delta^* \|\varphi\|_{L^2(\bd E)}.\]
\end{lemma}

We now recall the definition of inner and outer ball condition.
\begin{definition}\label{unifball}
    We say that a open set $E \subset \T^N$ satisfies a \textit{uniform inner (respectively outer) ball condition} with radius $r$ if there exists $r>0$ such that for every $x \in \bd E$ there exists a ball $B_r(y) \subset E$ (resp. $B_r(y) \subset E^c$) with $x \in \bd B_r(y)$.
\end{definition}
Note that all sets $E\subset \T^N$ of class $  C^{1,1}$ satisfy a uniform inner and outer ball condition (see e.g. \cite{Dal}).  Arguing as in the proof of \cite[Lemma 3.8]{AFM}, we can prove the following result.
\begin{lemma}\label{lemma1.3}
Let $E \subset \T^N$ be of class $C^{\infty}$ and let $m>0$. There exists $\eta=\eta(m,E)>0$ such that, for every $k\in \N$, $ u\in C^k(\bd E)$ with $\|u\|_{C^k(\bd E)}\le m$, $\|u\|_{C^0(\bd E)}\le \eta$ and for every $\sigma\in\T^N$ with $ |\sigma|\le \eta,$ then $E_u+\sigma$ can be written as a normal deformation of $E$ induced by a function $v:\bd E\to \bd E$ such that
\[\|v\|_{C^0(\bd E)}\le 2 \eta, \quad \|v\|_{C^k(\bd E)}\le C (\|u\|_{C^k(\bd E)}+|\sigma|), \]
where $C=C(E)>0$.
\end{lemma}
\begin{proof}
Being the set $E$ smooth, it satisfy the uniform inner and outer ball condition,  hence there exists a positive radius $r>0$ such that the signed distance $\sd_E$ from the set $E$, defined by
\begin{equation*}
\sd_E(x)= 
\begin{cases}
\dist_{\bd E}(x) &\text{if $x \in E^c$}\\
-\dist_{\bd E}(x) &\text{if $x \in E$,}
\end{cases} 
\end{equation*}
is a function of class $C^{\infty}$ (from the regularity of $\bd E$) in the $r$-tubular neighborhood $(\bd E)_{r}$, that is $ (\bd E)_r \coloneqq  \left\{ x \colon \, \dist_{\bd E}(x)<r\right\}$ (for further properties of the distance function see \cite[section 14.6]{gilbarg1977elliptic}). 
  Since, for some $k\ge 2$, $u$ has $C^k$-norm bounded by $m$, we also have $\|u\|_{C^{1,1}(\bd E)}\le m$. Then, there exists a radius $\rho=\rho(m,E)$ such that $\bd E_u$ satisfies a uniform inner and outer ball condition of radius $\rho.$
 We can assume without loss of generality that $\rho<r$. 
 
 We now let $\eta\le\rho/2$ to be chosen later, take any $|\sigma|< \eta$ and set $F=E_u+\sigma$. Clearly, $F$ still satisfies a uniform inner and outer ball condition of radius $\rho$. Then, for every $y\in\bd F$ there exists $x\in\bd E_u$ such that $y=x+\sigma$, hence we have
\begin{equation*}\label{neigh}
    \dist_{\bd E}(y)\le |\sigma|+\dist_{\bd E}(x)< \eta+\|u\|_{C^0(\bd E)}\le 2\eta,
\end{equation*}
and in particular $\bd F\subset (\bd E)_{2 \eta} \subset (\bd E)_r$.
We now define the map $T_u:\bd E\to \bd E$ as
\begin{equation}\label{def T}
T_u(x)\coloneqq \pi_{ E}(x+u(x)\nu_E(x)+\sigma)= y-\sd_E(y) \nabla \sd_E(y), \end{equation}
where $\pi_{E}$ is the projection map on $\bd E$ and $y=x+u(x)\nu_E(x)+\sigma \in \bd F$.  By choosing $\eta$ smaller, by interpolation, it holds $\|u\|_{C^1(\bd E)} + | \sigma | < \frac{1}{2} $,  which implies that the function $x \mapsto x+ u(x)\nu_{E}(x)+\sigma$ is a diffeomorphism (since it is a small perturbation of the identity). 
Moreover, since $E$ is of class $C^{\infty} $ (and possibly for $\eta $ smaller),   $\pi_{ E} \big|_{\bd F}:\bd F \to \bd E$ is a diffeomorphism of class $C^k$, $C^k$-close  to the identity. Therefore,  $T_u\in C^k(\bd E)$ and, by \eqref{def T}, we get
\begin{equation}\label{stimaT1}
     \|T_u-I\|_{C^k(\bd E)}\le C(\|u\|_{C^k(\bd E)}+|\sigma|).
\end{equation}
Moreover, using again \eqref{def T} and the invertibility of the map $x \mapsto x+u(x)\nu_E(x)+\sigma$, we obtain 
\begin{equation}\label{stimaT2}
      \|T^{-1}_u-I\|_{C^k(\bd E)}\le C(\|u\|_{C^k(\bd E)}+|\sigma|).
\end{equation}
Using the fact that $T_u$ is a diffeomorphism and \eqref{def T}, we can find a function $v:\bd E\to\R$ such that $F$ is the normal deformation of $E$ induced by $v$, more precisely for every $ x\in\bd E$ it holds
\begin{equation*}\label{rel}
x+u(x)\nu_E(x)+\sigma=T_u(x)+v(T_u(x))\nu_E(T_u(x)).
\end{equation*}
Finally, using the above expression and the bounds in \eqref{stimaT1} and \eqref{stimaT2}, we conclude that \[\|v\|_{C^k(\bd E)}\le \|T_u^{-1}\|_{C^k(\bd E)}(\|u\|_{C^k(\bd E)}+|\sigma|+\|T_u-I\|_{C^k(\bd E)})\le C(\|u\|_{C^k(\bd E)}+|\sigma|),\]
for some constant $C=C(E)>0$.
\end{proof}

Let $E, F\subset \T^N$ be measurable sets. We define a $L^1-$distance between $E,F$ modulo translations  (also known as the Fraenkel asymmetry of the set $E$ related to $F$) as
\begin{equation*}
    \alpha(E,F) \coloneqq  \min_{x \in \T^N} |E \triangle (F+x)|.
\end{equation*}

The following quantitative isoperimetric inequality has been proved in \cite{AFM}. As a consequence of this result, strictly stable sets are of class $C^{\infty}$ (see \cite{Mag}). 

\begin{theorem}[{\cite[Corollary 1.2]{AFM}}]
\label{coroll 1.2}
    Let $E \subset \T^N$ be a strictly stable set. Then, there exist $\eta= \eta(E)$, $C=C(E)>0$ such that 
    \[ C\alpha^2(E,F)\leq P(F)-P(E) \]
    for all $F\subset \T^N$ with $|F|=|E|$ and $\alpha(E,F)< \eta$.
\end{theorem}

We now recall the quantitative version of  Alexandrov's theorem proved in  \cite[Theorem~1.3]{DeKu}, which can be also seen as a Lojasiewicz-Simon inequality with sharp exponents. It will be the main tool to prove the exponential stability of the geometric flows considered. We slightly rephrase the conclusion as it will be more useful in the following.

\begin{theorem}[{\cite[Theorem 1.3]{DeKu}}]
  \label{teo alex}
Let $E \subset \T^N$ be a strictly stable critical set.
  There exist $\delta^* \in(0,1/2)$ and $C=C(E)>0$ with the following property: for any $f\in C^1(\bd E)\cap H^2(\bd E)$ such that $\|f\|_{C^1(\bd E)}\leq \delta^*$ and satisfying
  \begin{equation}
       |E_f|=|E|,\qquad \left |\int_{\bd E} f\nu_E\ud \mathcal H^{N-1}\right|\leq \delta^* \|f\|_{L^2(\bd E)},
       \label{ipotesi alex}
  \end{equation}
 setting $\mathscr{H}_{E_f}(x)=\HHH_{E_f}(x+f(x)\nu_E(x))$ for $x \in \bd E,$ we have
  \begin{equation}
      \|f\|_{H^1(\bd E)}\leq C\| \mathscr{H}_{E_f} - \bar{\mathscr{H}_{E_f}}
      \|_{L^2(\bd E)}.
      \label{e:aleq_gen}
  \end{equation}
\end{theorem}
\begin{remark}
      Note that equation \eqref{e:aleq_gen} in particular implies that, under the hypotheses of Theorem \ref{teo alex}, for any $\lambda\in\R$ it holds 
     \begin{equation}
      \|f\|_{H^1(\bd E)}\leq C\|\mathscr H_{E_f}-\lambda\|_{L^2(\bd E)}.
      \label{e:aleq_gen2}
  \end{equation}
\end{remark}
  
We conclude this section by recalling the Poincaré and Gagliardo-Nieremberg inequalities on smooth hypersurfaces (see \cite{Aub} for instance).

\begin{lemma}\label{Poincaré inequality}
    Let $\Sigma\subset \T^N$ be a smooth closed hypersurface  and $f\in H^1(\Sigma)$. There exists $C=C(\Sigma)>0$ such that
    \[ \| f- \bar f\|_{L^2(\Sigma)}\le C \| \n_\tau f \|_{H^1(\Sigma)}, \]
    {where we recall $\n_\tau f := \n f - (\n f \cdot \nu_\Sigma ) \nu_\Sigma$.}
\end{lemma}

\begin{theorem}
Let $\Sigma\subset \T^N$ be a smooth closed hypersurface. Let $l, \ m, \ k \in \N$ be such that $1\le l <m$, and let $1\leq r\le \infty$. There exists a constant $C$, depending on these constants and on $\Sigma$, with the following property: for every $u \in W^{l,p}( \Sigma)$ we have
$$\| \nabla^l u(\cdot,t)\|_{L^p(\Sigma )} \leq C \| u(\cdot,t) \|^{\theta}_{W^{m,r}(\Sigma)} \| u(\cdot,t) \|^{1-\theta}_{L^{q}(\Sigma)},  $$
where
\[\frac{1}{p}=\frac{l}{N-1}+\theta \left (\frac{1}{r}-\frac{m}{N-1} \right )+(1-\theta) \frac{1}{q}\]
for all $\theta \in [l/m,1)$ for which $p$ is nonnegative.
\end{theorem}

\subsection{Short-time existence for the mean curvature flow}

Given $T>0$ and $E_0 \subset \T^N$ an open smooth set, the \textit{volume-preserving mean curvature flow} in $[0,T)$ starting from $E_0$ is the family of sets $(E_t)_{0\le t < T}$ whose outer normal velocity is given by
\begin{equation}\label{defmcf}
V_t(x)= -\HHH_{E_t}(x) +\bar{\HHH}_{E_t}, \quad x \in \bd E_t, \ t \in (0,T).
\end{equation}
We remark that this equation should be intended as follows: there exist a smooth open set $E \subset \T^N$ and a $1$-parameter family of smooth diffeomorphism $\Phi_t:\bd E \to\T^N$ given by $\Phi_t(x)=x+u(x,t) \nu_E(x)$, such that $\Phi_0(\bd E)= \bd E_0$,  $\Phi_t(\bd E)=\bd E_t$, and
\begin{equation*}
\partial_t u(x,t) \nu(x) \cdot \nu_{E_t}(\Phi_t(x))= -\HHH_{E_t}(\Phi_t(x)) + \bar{\HHH}_{E_t}, \quad x \in \bd E, \ t \in (0,T).
\end{equation*}

Assuming that the flow starting from $E_0$ exists, following classical computations (see for instance \cite{Manlib}) one can deduce that the evolution equation satisfied by $u$ is 
\begin{equation*}
    \begin{split}
        \partial_t u = \Delta_{ E} u +  \langle A(x,u,\nabla u), \nabla^2 u \rangle +J(x,u,\nabla u)+\HHH_E, 
    \end{split}
\end{equation*}
where $\Delta_{E}$ is the Laplace-Beltrami operator on ${\bd E}$, $A$ is a smooth tensor such that $A(\cdot,0,0)=0$,  and   $J$ is a smooth function.

In order to prove the stability of such flow, we need the following short-time existence result. 

\begin{theorem}\label{teo short time MCF}
    Let $\varepsilon>0$, let $\beta\in (0,1)$ and let $E \subset \T^N$ be a smooth open set.
    There exists $\delta=\delta(\varepsilon,E,\beta)>0$ with the following property: if $E_0$ is the normal deformation of $E$ induced by $u_0 \in C^{1,1}(\bd E)$, $\|u_0\|_{C^{1,1}(\bd E)}\le \delta$, and $|E_0|=|E|$, then there exists $T>0$, which only depends on $E$, $\beta$ and the bound on $\|u_0\|_{C^{1,1}(\bd E)}$, such that the volume preserving mean curvature flow $E_t$ starting from $E_0$ exists in $[0,T)$, the sets $E_t$ are  normal deformation of $E$ induced by $u(\cdot,t) \in C^{\infty}(\bd E)$ for all $t \in (0,T)$, and 
    \begin{equation}\label{stima1}
        \sup_{t\in (0,T)} \|u(\cdot,t)\|_{C^{1,\beta}(\bd E)} \le \varepsilon.
    \end{equation}
    Moreover, for every $k \in \N $, there exist two constants  $c_k=c_k(N)>0$ and $C_k=C_k(E)>0$ such that
    \begin{equation}\label{stima2}
         \sup_{t\in (0,T)} t^{c_k}\|\n^{k+2} u(\cdot,t) \|_{C^0(\bd E)}\le C_{k} (\| u_0 \|_{C^{1,1}(\bd E)}+1). 
    \end{equation}
\end{theorem}

We remark that the proof of this result  is classical and can be derived from the Schauder estimates for quasi-linear parabolic equations, as  $u$ solves a lower-order, nonlinear perturbation of the heat equation. In the following section we will provide a brief outline of the proof for an analogous short-time existence result for the surface diffusion flow (see Theorem~\ref{shortsurf}). Similar
and simplified arguments would prove the previous result for the mean curvature flow, which is a second order flow.

For the sake of completeness, we provide here an alternative proof of Theorem~\ref{teo short time MCF} which follows from some results found in the literature. Even if these results are  shown in the ambient space $\R^N$, the same arguments can be repeated in the flat torus.  The first part of the Theorem is the short-time existence result of \cite{ES}.

\begin{theorem}[{\cite[Main Theorem]{ES}}]\label{teo short time ES}
    Let $E\subset \T^N$ be a smooth open set and $\beta\in (0,1)$. There exists $\delta=\delta(E,\beta)>0$  with the following property: if $E_0$ is the normal deformation of $E$ induced by $u_0 \in C^{1,1}(\bd E)$, $\|u_0\|_{C^{1,1}(\bd E)}\le \delta$, and $|E_0|=|E|$, then there exists $T>0$, only depending on $E$, $\beta$ and the bound on $\|u_0\|_{C^{1,1}(\bd E)}$, such that the volume-preserving mean curvature flow $E_t$ starting from $E_0$ exists in $[0,T)$, and the sets $E_t$ are normal deformations induced by $u(\cdot,t) \in C^{\infty}(\bd E)$ for all $t \in (0,T)$. Furthermore, the mapping $(t,E_0)\mapsto E_t$ is a local smooth semiflow on $C^{1,\beta}(E)$.
\end{theorem}

We remark that the local smooth semiflow property in particular implies that $\|u(\cdot)\|_{C^{1,\beta}}$ depends continuously on $\|u_0\|_{C^{1,\beta}}$ (see for instance \cite[pag. 66]{Amalib}). In particular, for every $\e>0$ there exists $\delta(E,\e,\beta)>0$ and $T(E,\e,\beta)>0$ such that if $\|u_0\|_{C^{1,\beta}}\le \delta$ then
\beq
    \|u(\cdot,t)\|_{C^{1,\beta}}\le \e \quad \text{for every } t\in (0,T).
\eeq

In order to obtain the higher-order regularity inequalities, we apply some curvature estimates obtained recently in \cite{JN}.

\begin{theorem}[{\cite[Theorem 1.1]{JN}}]\label{JN main}
Assume that $E_0\subset \R^N$ is an open bounded set satisfying a uniform inner and outer ball condition with radius $r$. Then, there exists a time $T=T(r,N)>0$ such that the volume preserving mean curvature flow $E_t$ starting from $E_0$ exists in $[0,T)$ and it satisfies a uniform inner and outer ball condition of radius $r/2$. Moreover, it is smooth in $(0,T)$ and satisfies for every $k\in \N$
\begin{equation}\label{stime Jul}
\sup_{t\in(0,T)} \left( t^k\|\HHH_{E_t}\|_{H^k(\bd E_t)}^2 \right)\le C_k,  
\end{equation}
where $C_k$ depends on $k,|E_0|,r$. 
\end{theorem}

Before proving the short time existence result, we  remark a classical result concerning the uniform ball condition.
 \begin{remark}\label{rmk ball C11}
    Let $E$ be a smooth set satisfying a uniform ball condition of radius $r_E$. Then every small $C^{1,1}$-normal deformations of $E$ satisfy a uniform ball condition of radius $r\approx r_E$. Indeed, it is easy to see that if $E_f$ is the normal deformation of $E$ induced by $f\in C^{1,1}(\bd E)$, then the Hausdorff distance between $E$ and $E_f$ is bounded by $\|f\|_{C^0(\bd E)}.$ Furthermore, since $\n \sd_{E_f}=\nu_{E_f}$ and $\nu_{E_f}$ can be written  as  
    \begin{equation}
	\nu_{E_f}=\left(\nu_E-\sum_{i=1}^{N-1} \frac{\n f\cdot v_i}{1+\kappa_i f}v_i\right)  \left(  1+ \sum_{i=1}^{N-1} \dfrac{(\n f\cdot v_i)^2}{(1+\kappa_i f)^2} \right)^{-1/2},
	\label{vers normale E_f}
\end{equation}
    where the family $\{v_i\}_{i=1,\dots,N-1}$ denotes an orthonormal frame of the tangent space on $\bd E$ (see \cite[eq. (3.3)]{DeKu}),  by  differentiating \eqref{vers normale E_f} one can see that 
    \[  \| \sd_{E_f}-\sd_E \|_{C^{1,1}(\bd E)}\le C_E \|f\|_{C^{1,1}(\bd E)}, \]
    which then implies that $E_f\to E$ in $C^{1,1}$ if $\|f\|_{C^{1,1}}\to 0$. Therefore, by   \cite[Theorem 2.6]{Dal} and \cite[Remark 2.7]{Dal} one infers that the radius $r$ of the uniform ball condition of the set $E_f$ depends continuously on $\|f\|_{C^{1,1}}$ when it is small enough. In particular, for every $\e>0$ there exists $\delta(r_E,\e)>0$ such that, if $\|f\|_{C^{1,1}}\le \delta$ then
    \beq
    |r_E-r|\le \e.
    \eeq
 \end{remark}

\begin{proof}[{Proof of Theorem~\ref{teo short time MCF}}]
By Theorem~\ref{teo short time ES} there exist a time $T'>0$ and a family of evolving functions $u(\cdot,t)$, which are smooth in $(0,T')$ and satisfy the inequality \eqref{stima1}. The second bound follows from classic elliptic regularity arguments that we now sketch. 
\\Fix $t\in (0,T')$, from the bound on $\sup_{t \in (0,T')}\| u \|_{C^{1,\beta}(\bd E)}$ and (up to rotations) for any given point $x=(x',x_{N})\in\bd E$ we can parametrize  in a cylinder $ C=B'_r(x)\times(-L,L)$   both $\bd E$ and $\bd E_t$ as graphs of smooth functions $g,g_t$. From Theorem~\ref{JN main} there exists a time $T''$ (depending on $E,\delta$ by Remark\eqref{rmk ball C11}) such that the evolving sets $E_t$ satisfy a uniform inner and outer ball condition of radius $r/2$ for any $t \in (0, T'')$.  Let us set $T = \min \{T', T''\}. $ From estimate \eqref{stime Jul} we get that
\[ \HHH_{E_t}=\div \left( \dfrac{\nabla g_t}{\sqrt{1+|\nabla g_t|^2}}\right) = \dfrac{1}{\sqrt{1+|\nabla g_t|}}\left( I-\dfrac{\nabla g_t\otimes \nabla g_t}{1+|\nabla g_t|^2}  \right):\nabla^2  g_t \]
is bounded in $L^{2}(B'_r(x'))$ by a constant which depends on $\vert E_0 \vert ,T,r$. 
Then, by uniform geometric Calderon-Zygmund inequality (see \cite[Section 3]{DDM} or \cite[Lemma 7.2]{AFM}) we deduce that, for some $\rho<r$, in the ball $B'_{\rho}(x')$ the function $g_t$ is bounded in $H^{2}(B'_{\rho}(x'))$ by a constant,
depending only on the $L^2$-bound on $\HHH_{E_t}$, the norm of the coefficients of the elliptic operator, which are in turn bounded by $\|u_0\|_{C^{1,1}}$ thanks to the previous step. Iterating this procedure, we bound the higher norms $H^k(B'_{\rho}(x'))$ of $g_t$, for every $k\in\N$. 
Then, we conclude by means of Sobolev embeddings and by a covering argument.
 \end{proof}

\subsection{Short-time existence for the surface diffusion flow}

We now consider the evolution called \textit{surface diffusion flow}, defined by
\begin{equation}\label{defsdf}
V_t(x)=\Delta_{E_t} \HHH_{E_t}(x), \quad x \in \bd E_t, \ t \in (0,T).
\end{equation}
As for the mean curvature flow, the equation above means that there exist a smooth open set $ E \subset \T^N$ and a $1$-parameter family of smooth diffeomorphism $\Phi_t:E \to \T^N$ such that $\Phi_t(x)=x+u(x,t) \nu_E(x)$, $\Phi_t(\bd E)=\bd E_t$ and
\begin{equation*}
\partial_t u(x,t) \nu_E(x) \cdot \nu_{E_t}(\Phi_t(x))= \Delta_{E_t}\HHH_{E_t}(\Phi_t(x)).
\end{equation*}
Assuming that the diffeomorphisms above exist, arguing as in \cite[pag. 21]{Manlib}, one can deduce that the evolution equation satisfied by $u$ is 
\begin{equation}\label{eqqqqqq}
\begin{split}
    \partial_t u =& -  \Delta^2_{E_t} u - \dfrac 1{\nu_E\cdot \nu_{E_t}} \Delta_{E_t} (\nu_E\cdot \nu_{E_t} ) \Delta_{E_t} u + \dfrac 1{\nu_E\cdot \nu_{E_t}}\Delta_{E_t} P(x,u,\n u)\\
    &= -  \Delta^2_{E_t} u + \tilde J(x,u,\n u,\n^2 u,\n^3 u),
    \end{split}
\end{equation}
where $P$ is a smooth function (assuming that $u$ and $\n u$ are small), the function $\tilde J$ can be written as 
\[ \tilde J(x,u,\n u,\n^2 u, \n^3 u)=  \la \tilde B_1,  \n^2 u \ra + \la \tilde B_2, \n^2 u\otimes \n^2 u \ra + \la \tilde B_3, \n^3 u \ra +\tilde b_4  \]
and $\tilde B_1,\tilde B_2,\tilde B_3$ and $\tilde b_4$ are tensor-valued, respectively scalar-valued functions depending on $(x,u,\n u)$ and smooth if their arguments are small enough. {Here, with a little abuse of notation, $\n$ denotes} the covariant derivative on $\bd E$.

On the other hand, linearizing the Laplace-Beltrami operator yields the evolution equation (compare with \cite[Section~3.1]{FusJulMor3D})
\begin{equation}\label{eq riscritta}
    \begin{split}
        \partial_t u = -\Delta_E^2 u + \langle A(x,u,\nabla u), \nabla^4 u \rangle + J(x,u,\nabla u,\nabla^2 u, \nabla^3 u), 
    \end{split}
\end{equation}
where  $A$ is a smooth $4$th-order tensor, vanishing when both $h$ and $\nabla h$ vanish, and $J$ is given by
\begin{equation}\label{formulaJestesa}
    \begin{split}
        J=&\langle B_1, \nabla^3 u \otimes \nabla^2 u \rangle + \langle B_2, \nabla^3 u\rangle + \langle B_3, \nabla^2 u \otimes \nabla^2 u \otimes \nabla^2 u \rangle \\
        &+ \langle B_4, \nabla^2 u \otimes \nabla^2 u \rangle + \langle B_5, \nabla^2 u \rangle + b_6,
    \end{split}
\end{equation}
where $B_i$, $i=1,\ldots 5$ and $b_6$ are smooth tensor-valued, respectively scalar-valued functions depending on $(x,u,\n u)$.

In this subsection we want to prove a short-time existence result for the surface diffusion flow, in particular we will obtain a priori estimates that will be used to prove the stability of the flow. We will follow the classical approach of linearization and fixed point to solve the nonlinear evolution problem, and then employ Shauder-type estimates to show higher order regularity of the flow. We will follow closely what has been done in \cite{FusJulMor3D}, combining it with the results of \cite{HZ}.

To start we recall some classical results concerning the Cauchy problem for the biharmonic heat equation on a smooth Riemannian manifold $\Sigma$ with metric $g$, which is the solution to the following problem
\begin{equation}\label{cauchy pb biharm}
    \begin{cases}
        \bd_t u=-\Delta_\Sigma^2 u+f(x,t)\quad&\text{on }\Sigma\times[0,\infty)\\
        u(\cdot,0)=u_0\quad&\text{on }\Sigma,
    \end{cases}
\end{equation}
once the functions $f,u_0$ are assigned.

\begin{theorem}[p.~251, {\cite[Theorem 2]{Fri}}]\label{thm301}
Given $( \Sigma, g)$ a smooth Riemannian manifold, there exists a unique biharmonic heat kernel with respect to $g$ denoted as $b_g \in C^\infty\big(\Sigma\times \Sigma \times (0, \infty)\big)$. Moreover let $T>0$, for any integers $k, p, q \geq 0 $ and for any $(x,y, t) \in \Sigma \times \Sigma \times(0,T)$ we have 
\begin{equation}\label{eqn32}
|\partial_t^k \nabla_x^{p} \nabla_y^{q} b_g(x, y,  t)|_g \leq C t^{- \frac{n + 4k + p+ q}{4}} \exp\{- \delta \big(t^{-\frac{1}{4}} d_g(x, y)\big)^{\frac{4}{3}} \},  
\end{equation}
where $|\cdot |_g=\sqrt{g(\cdot,\cdot)},$ $\nabla_x$ and $\nabla_y$ are covariant derivatives with respect to $g$, and the constants $C, \delta >0$ depend on $T$, $g$ and $p + q +4k$.  
\end{theorem}

Given the biharmonic heat kernel $b_g \in C^\infty\big(\Sigma \times \Sigma \times (0,\infty)\big)$ on $(\Sigma,g)$ and a function $u_0 \in C^0(\Sigma)$, we define for $(x, t) \in \Sigma \times (0,\infty)$
\begin{equation}\label{eqn4001}
S u_0 (x, t) = \int_\Sigma b_g(x, y, t) u_0(y) \mathrm d V_g(y) 
\end{equation} 
where $  V_g$ is the Riemannian volume form.
Hence, as usual, $Su_0$ is the solution to the homogeneous Cauchy problem
\begin{equation}\label{eqn4002}
\begin{cases}
      \partial_t v+ \Delta_\Sigma^2   v  = 0 \quad &\text{on }\Sigma\times (0,+\infty) \\
    v(\cdot,0)=u_0(\cdot)\quad \text{on }\Sigma.
\end{cases}
\end{equation} 
Moreover, since the biharmonic heat kernel is smooth for every $t>0$, we get  $S u_0 \in C^{\infty}\big(\Sigma \times (0,+\infty)\big)$.
We now collect some results, which are shown in \cite{HZ}, about the solution of \eqref{cauchy pb biharm}.  The following Schauder-type estimates on the solution of the homogeneous problem \eqref{eqn4002} can then be proved, see \cite[Theorem 3.8]{HZ}. In particular, we modify slightly the formulation of the result, to fit our purposes. One can inspect the proof of \cite[Theorem 3.8]{HZ} (see pag. 7487,7489 in particular) to check the result.

\begin{theorem}\label{thm001}
Suppose $u_0 \in C^{1,1}(\Sigma)$ and fix $T>0$. Then there exists $C_1(\Sigma,T)>0$ such that
\begin{equation}\label{schauder Su}
    \sup_{t\in (0,T)} \| | S u_0 |_g \|_{C^{1,1}(\Sigma)}\le C_1 \|u_0\|_{C^{1,1}(\Sigma)},
\end{equation}
Furthermore, for any $l, k \in \N$, we have 
\begin{equation}\label{higher order estimates S u}
\sup_{t \in (0,T)}t^{l + \frac{k}{4}} \bigg \| \Big | \left (\bd_t \right )^l \nabla_g^{k+2} S u_0(t) \Big |_g \bigg\|_{C^0(\Sigma)} \leq C_{l,k} \|u_0\|_{C^{1,1}(\Sigma)}, 
\end{equation}
for some constants $C_{l,k} >0$ depending on $l$, $k$, $\Sigma$ and $T$.
\end{theorem}

In order to study the evolution problem \eqref{eq riscritta} we introduce the following two Banach spaces. Fix $0< T < \infty$ and $0< \beta < 1$. We define 
\begin{equation}\label{space1}
Y_T \coloneqq  \left \{u \in C^{0}\big(\Sigma \times (0,T) \big) : \|u\|_{Y_T} < \infty \right \},
\end{equation}
where
\begin{equation}
\begin{split}
\|u\|_{Y_T} \coloneqq &  \sup_{t\in(0,T) } \left ( t^{\frac{1}{2}} \|u(\cdot,t)\|_{C^0(\Sigma)} + t^{\frac{1}{2} + \frac{\beta}{4}} [u(\cdot,t)]_{C^{\beta}(\Sigma)} \right) \\
&+ \sup_{(x,t) \in \Sigma \times (0,T)} \sup_{0< h < T-t} t^{\frac{1}{2} + \frac{\beta}{4}}\frac{|u(x, t+h) - u(x,t)|}{|h|^{\frac{\beta}{4}}}
\end{split}
\end{equation} 
and $[\cdot]_{C^\beta}$ is the usual H\"older seminorm. Similarly, we introduce the space 
\begin{equation}\label{space2}
X_T \coloneqq  \left \{u \in C^0(\Sigma \times (0,T)) : u(\cdot,t) \in C^4(\Sigma), \ \|u\|_{X_T} < \infty \right \},
\end{equation}
where 
\begin{equation}\label{norma X_T}
\begin{split}
\|u\|_{X_T} \coloneqq & \sup_{t\in (0,T)} \Big( \sum_{k=0}^4  t^{- \frac{1}{2} + \frac{k}{4}} \|\nabla^k  u(\cdot,t)\|_{C^0(\Sigma)}+  t^{\frac{1}{2}+ \frac{\beta}{4}} [\nabla^4 u(\cdot,t)]_{C^{\beta}(\Sigma)} \\
& + t^{\frac{1}{2}} \|\bd_t  u(\cdot,t)\|_{C^{0}(\Sigma)} +  t^{\frac{1}{2}+\frac{\beta}{4}} [\bd_t  u(\cdot,t)]_{C^{\beta}(\Sigma)} \Big) \\
&+ \sup_{(x,t) \in \Sigma \times (0,T)} \sup_{0< h < T-t} t^{\frac{1}{2}+\frac{\beta}{4}}\frac{|\nabla^4u (x, t+h)- \nabla^4 u(x,t)|_g}{|h|^{\frac{\beta}{4}}}\\
& + \sup_{(x,t) \in \Sigma \times (0,T)} \sup_{0< h < T-t} t^{\frac{1}{2}+\frac{\beta}{4}}\frac{|\bd_t  u (x, t+h)- \bd_t  u (x,t)|}{|h|^{\frac{\beta}{4}}}. 
 \end{split}
\end{equation}

\begin{proposition} The spaces $(Y_T, \Vert \cdot \Vert _{Y_T})$ and $(X_T, \Vert \cdot \Vert _{X_T})$ are Banach spaces.
\end{proposition}
The proof of the completeness of the spaces $Y_T$ and $X_T$ is standard, indeed one can prove directly that all Cauchy sequence converge to a function in the space and the candidate limit is obtained using a diagonal argument. 
\begin{remark} \label{rmk norms}
Since the norm $\sum_{k=0}^4\|\n^k u\|_{C^0}$ is equivalent to the norm $\|u\|_{C^0}+\|\n^4u\|_{C^0}$ for $C^4(\Sigma)$, we have that the norm $\|\cdot\|_{X_T}$ defined in \eqref{norma X_T} is equivalent to the following norm 
\begin{equation*}
\begin{split}
\|u\|_{X_T}^{'} : = & \|u\|_{X_T} + \sum_{k= 0}^3 \sup_{(x,t) \in \Sigma \times (0,T)} \sup_{0< h < T-t} t^{- \frac{1}{2} + \frac{k}{4}+\frac{\beta}{4}}\frac{|\nabla^k u (x, t+h)- \nabla^k u(x,t)|_g}{|h|^{\frac{\beta}{4}}}.
\end{split}
\end{equation*}
\end{remark}

Now we study the nonhomogeneous initial value problem
\begin{equation}\label{eqn7001}
\begin{cases}
\bd_t u  + \Delta_\Sigma^2 u = f \quad &\text{ on } \Sigma \times (0,T) \\
u(\cdot, 0) = 0 &\text{ on }  \Sigma,
\end{cases}
\end{equation}
where $f $ is a function on $\Sigma \times (0,T)$. Given the biharmonic heat kernel $b_g \in C^\infty\big(\Sigma \times \Sigma \times (0,T)\big)$ on $(\Sigma,g)$, the solution (if it exists) to the nonhomogeneous problem \eqref{eqn7001} should be given by Duhamel's principle
\begin{equation}\label{duhamel}
V f (x, t) \coloneqq  \int_{0}^t \int_\Sigma b_g(x, y, t-s) f(y, s) \mathrm d V_g(y) \mathrm d s, 
\end{equation}
and, for every $\lambda>0$, $Vf \in C^{\infty}(\Sigma \times (\frac{\lambda}{2},\lambda)).$

We then recall the following fundamental Schauder-type estimates proved in \cite{HZ} on solutions of \eqref{eqn7001} (see \cite[Remark~3.12]{HZ} for the final comments on the constant $C$).

{\begin{theorem}[{\cite[Theorem~3.10]{HZ}}]\label{thm701}
Fix $0 < T < \infty$,
if $f \in Y_T$, then $Vf \in X_T$ and there exists a constant $C>0$ depending on $\Sigma,T$ such that
\begin{equation}
\|V f\|_{X_T} \leq C \|f\|_{Y_T}. 
\end{equation}
Moreover, equation $(\partial_t  + \Delta_\Sigma^2 ) Vf = f$ holds in the classical sense on $\Sigma \times (0,T)$ and thus $V f\in C^\infty(\Sigma\times (0,T))$. 
\end{theorem}}

We now turn our attention to the evolution equation \eqref{eq riscritta}, and use the results above for the particular choice $\Sigma= \bd E$ with the Riemaniann metric induced by the Euclidean one. We consider the map 
\beq\label{error term}
    f[u](x)\coloneqq \langle A(x,u,\nabla u), \nabla^4 u \rangle + J(x,u,\nabla u,\nabla^2 u, \nabla^3 u),
\eeq
where $A,$ $J$ are the operators defined in \eqref{eq riscritta}.
We now provide the fundamentals estimates on $f[u]$, which represents the nonlinear error generated linearizing \eqref{eq riscritta}.

\begin{lemma}\label{lemma est Jh}
For any $\e,\,m>0$ there exist  $T, \,\delta>0$ depending on $E,\e$ with the following properties. For every $u_0\in C^{1,1}(\Sigma)$ and $\psi \in X_T$ satisfying $\|\psi\|_{X_T}\le m$ it holds
\beq\label{in space}
    f[\psi+S u_0]\in Y_T. 
\eeq
Moreover, if $\| u_0\|_{C^{1,1}(\Sigma)}\le \delta$ it holds
\begin{equation}\label{est Jh}
    \|f[S u_0]\|_{Y_T}\le \e  ( \| u_0 \|_{C^{1,1}(\Sigma)}+ 1).
\end{equation}
Finally,    $\psi_1,  \psi_2 \in X_{T}$  satisfying   $ \|\psi_i\|_{X_{T}}\le m$,  it holds
\begin{equation}\label{contrazione}
 \|f[\psi_1+S u_0]-f[\psi_2 + S u_0]\|_{Y_{T}} \leq \e \|\psi_1-\psi_2 \|_{X_{T}}.
\end{equation} 
\end{lemma}

\begin{proof}
Let $T<1$ to be chosen later and fic $\e,m>0$.
We prove only equation \eqref{est Jh}, giving a sketch of the proof for \eqref{contrazione} and \eqref{in space} as they are analogous; we also drop the dependence on the set $E$ in the norms.  For clarity of exposition, we prove the results for the simplified error term 
\beq\label{sempli err}
    \tilde f[u](x,t)\coloneqq \la A(x,u(x,t),\n u(x,t)), \n^4 u(x,t)  \ra +\la B, \n^3 u(x,t)\otimes \n^2 u(x,t) \ra,
\eeq
where $B$ is a (constant) tensor of the same dimension of $\n^3 u\otimes \n^2 u$ with $\| B \|< 1$.  The general case is explained in the appendix, but follows by analogous computations. 
We will also write $A(x,t)$ and assume implicitly the dependence on $u,\n u$. 

Firstly, we prove \eqref{est Jh}. In what follows we use the short-hand notation $u=S u_0$. 
From the definition of $\tilde f[\cdot]$ we have
\beq 
\begin{split} \label{equazione qua}
    \| \tilde f[u] \|_{C^0}&\le \|A\|_{C^0}\|\n^4 u\|_{C^0}+ \|\n^3 u\|_{C^0}\|\n^2 u\|_{C^0},\\
     [\tilde f[u]]_{C^\beta}&\le  \|\n^4 u\|_{C^0} \sup_{\tau\in\T^N} \left( |\tau|^{-\beta}|A(x+\tau,t)-A(x,t)|\right) + \|A\|_{C^0} [\n^4 u]_{C^\beta}\\
     &\  + [\n^3 u]_{C^\beta}\|\n^2 u\|_{C^0}  +  \|\n^3 u\|_{C^0}[\n^2 u]_{C^\beta}.
\end{split}
\eeq
Then, we multiply by $t^{\frac 12}$ the first equation in \eqref{equazione qua} to get
\[  t^{\frac12}\| \tilde f[u] \|_{C^0}\le \|A\|_{C^0}  t^{\frac12} \|\n^4 u\|_{C^0}+  t^{\frac14} t^{\frac14}\|\n^3 u\|_{C^0}  \|\n^2 u\|_{C^0}. \]
By \eqref{higher order estimates S u}, with the choice of $l=0$, $k=0,1,2$, we have that all the terms $t^{\frac12} \|\n^4 u\|_{C^0}$, $ t^{\frac14}\|\n^3 u\|_{C^0}$ and $\|\n^2 u\|_{C^0}$ are bounded by $\|u\|_{C^{1,1}}$ {(times a constant that depends on $E$ which we can suppose equal to one for simplicity)}. We now fix $\delta>0$ sufficiently small, depending on $\varepsilon$ and $E$, so that $\|A\|_{C^0}$ is bounded by $\varepsilon$, which can be done since $A$ is a smooth tensor and $A(\cdot,0,0)=0$. Finally, taking $T$ small enough, depending on $\varepsilon$ and $E$, we conclude 
\[ \sup_{t\in (0,T)} t^{\frac12}\| \tilde f[u] \|_{C^0}\le  \e \|u_0\|_{C^{1,1}}. \]
Therefore, taking into account the full expression for the error term $f[u]$  given by \eqref{error term}, one can show that
\[ \sup_{t\in (0,T)} t^{\frac12}\|  f[u] \|_{C^0}\le C\e \left( \|u_0\|_{C^{1,1}}+1\right), \]
where the last constant comes from the term $b_6$.

Concerning the H\"older seminorm in space, we first remark that
    \[   \sup_{\tau\in \T^N} \frac{|A(x+\tau,t)-A(x,t)|}{|\tau|^{\beta}}\le [A(\cdot, u,\n u)]_{C^\beta}+ \| \bd_2 A \|_{C^0}[u]_{C^\beta}+\| \bd_3 A \|_{C^0}[\n u]_{C^\beta}, \]
  where  $ \partial_2 A$ and $\partial_3 A$ denote the derivative of $A(x,y,z)$ with respect to the second and third components.
 Therefore, employing again the bounds in   \eqref{schauder Su} and \eqref{higher order estimates S u} we can bound
    \beq
        t^{\frac 12}\|\n^4 u\|_{C^0} \sup_\tau\frac{|A(x+\tau,t)-A(x,t)|}{|\tau|^{\beta}}\le \varepsilon \| u_0\|_{C^{1,1}},
    \eeq
    where we took $\delta>0 $ sufficiently small, depending on $\varepsilon$ and $E$, such that 
    \begin{equation*}
        [A(\cdot, u,\n u)]_{C^\beta}+ \| \bd_2 A \|_{C^0}[u]_{C^\beta}+\| \bd_3 A \|_{C^0}[\n u]_{C^\beta} \leq \varepsilon,
    \end{equation*}
    which is possible since $A$ is smooth and  $ A(\cdot,0,0)=0$.
    Thus, multiplying by $t^{\frac 12+\frac \beta4}$ the second equation in \eqref{equazione qua} we obtain
    \begin{equation}
        \begin{split}\label{eq above}
             t^{\frac12+\frac \beta 4}[\tilde f[u]]_{C^\beta}&\le t^{\frac \beta 4}\varepsilon\| u_0\|_{C^{1,1}} + \|A\|_{C^0}t^{\frac12+\frac \beta 4} [\n^4 u]_{C^\beta}\\
        &\ +t^{\frac 14} t^{\frac14+\frac \beta 4} \|\n^3 u\|_{C^{\beta}}\|\n^2 u\|_{C^0} +  t^{\frac14} t^{\frac 14}\|\n^3 u\|_{C^{0}} t^{\frac \beta 4} \|\n^2 u\|_{C^\beta}.
        \end{split}
    \end{equation}
    Then, all the terms in \eqref{eq above} with the norms of $u$ can be bounded employing \eqref{schauder Su} and \eqref{higher order estimates S u}, thus we can make the right-hand side above as small as needed taking   $T,\delta$ small enough. Analogous calculations show a similar inequality for the complete error term $f[u].$

    Finally, we show how to bound  the H\"older seminorm in time appearing in $\|\tilde f[u]\|_{Y_T}$. We fix $t\in (0,T),h\in (0,T-t)$. To ease notation, we omit to write the evaluation at $x$ in the following. We have by the very definition of $ \tilde f[u](t)$ that 
    \begin{equation*}
     \begin{split}
          | \tilde f[u]&(t+h)- \tilde f[u](t)| \\
        & \le | \la A( u( t+h),\n u( t+h)), \n^4 u( t+h)  \ra - \la A( u( t),\n u( t)), \n^4 u( t)  \ra | \\
        & \; + | \la B,\left( \n^3 u( t+h)\otimes \n^2 u( t+h) \right)\ra - \la B, \left( \n^3 u( t)\otimes \n^2 u( t) \right)\ra |.
     \end{split}   
    \end{equation*}
     Now by the triangular inequality we obtain \begin{equation}\label{triangologranginequality1}
        \begin{split}
           & | \la A(u( t+h),\n u( t+h)), \n^4 u( t+h)  \ra - \la A( u( t),\n u( t)), \n^4 u( t)  \ra | \\
             & \leq \|A \|_{C^0} | \n^4 u(t+h)- \n^4 u(t) | + \| \partial_3 A \|_{C^0} | \nabla u(t+h)- \nabla u(t) | \| \nabla^4 u(t)\|_{C^0} \\
             & +\| \partial_2 A \|_{C^0} | u(t+h)-u(t) | \| \nabla^4 u \|_{C^0},
        \end{split}
    \end{equation} 
     and analogously
\begin{equation}\label{triangologranginequality2}
        \begin{split}
           &| \la B,\left( \n^3 u( t+h)\otimes \n^2 u( t+h) \right)\ra - \la B,\left( \n^3 u(x,t)\otimes \n^2 u(x,t) \right)\ra | \\
           & \leq  |\n^3 u(t+h)-\n^3 u(t)|\|\n^2 u\|_{C^0}
         +\|\n^3 u \|_{C^0} |\n^2 u(t+h)-\n^2 u(t)|.
        \end{split}
    \end{equation}
    Therefore from formulas \eqref{triangologranginequality1} and \eqref{triangologranginequality2},  we obtain 
    \begin{align*}
        |\tilde f[u](t+h)&- \tilde f[u](t)| \\
        &\le  \left( \| \bd_2 A \|_{C^0} |u(t+h)-u(t)|+\|\bd_3 A\|_{C^0} | \n u(t+h)-\n u(t)|\right) \|\n^4 u(t)\|_{C^0}\\
        &\ + \|A\|_{C^0} |\n^4 u(t+h)-\n^4 u(t)|+ |\n^3 u(t+h)-\n^3 u(t)|\|\n^2 u\|_{C^0}\\
        &\ +\|\n^3 u \|_{C^0} |\n^2 u(t+h)-\n^2 u(t)|.
    \end{align*}
    Applying again \eqref{schauder Su}, \eqref{higher order estimates S u}, and using the smallness of $\|A\|_{C^0}$, we conclude \eqref{est Jh} by taking $T,\ \delta$ small enough. 
    
    Following the computations above one can easily prove that if $u_0\in C^{1,1}(\Sigma)$ and $\|\psi\|_{X_T}\le m$, it holds 
    \[ f[\psi+ Su_0]\in Y_T. \]
    The only difference is that, in addition to \eqref{schauder Su}, \eqref{higher order estimates S u} one can directly exploit the definition of $\|\cdot\|_{X_T}$ to obtain the required bounds. Also the proof for \eqref{contrazione} is essentially the same, only much more tedious to write. We show the computations only for the term $\sup_{t\in (0,T)} t^{1/2}\| \cdot \|_{C^0} $ appearing in the norm of $Y_T$ and for the simplified error term \eqref{sempli err}. For $u_i\coloneqq\psi_i+Su_0$ we can write
    \begin{align*}
        &\lvert  \tilde f[u_1] - \tilde f[u_2]  \rvert\\
        &=\left\lvert \la A(x,u_1,\n u_1), \n^4 u_1  \ra - \la A(x,u_2,\n u_2), \n^4 u_2  \ra +\la B, (\n^3 u_1\otimes \n^2 u_1 - \n^3 u_2\otimes \n^2 u_2)\ra    \right\rvert \\
        &\le\|\n^4 u_1\|_{C^0} \left (  \|\bd_1 A\|_{C^0} |  \psi_1-  \psi_2| +  \|\bd_2 A\|_{C^0} |\n \psi_1-\n \psi_2| \right)+\|A\|_{C^0} |\n^2 \psi_1-\n^2 \psi_2| \\
        &\ + \|\n^3 u_1 \|_{C^0} |\n^2 \psi_1-\n^2 \psi_2|+ \|\n^2 u_2 \|_{C^0} |\n^3 \psi_1-\n^3 \psi_2|.
    \end{align*}
    Multiplying the inequality above by $t^{\frac 12}$ we have
    \begin{align*}
        &t^{\frac 12}\lvert \tilde f[u_1] - \tilde  f[u_2]  \rvert \\
        &\le  \Big( \|\n^4 u_1\|_{C^0} \left ( t \|\bd_1 A\|_{C^0} +  t^{\frac 34} \|\bd_2 A\|_{C^0}  \right)+ t^{\frac 12 } \left (\|A\|_{C^0}+ \|\n^3 u_1 \|_{C^0}\right)\\
        &\ + t^{\frac 14 }\|\n^2 u_2 \|_{C^0}  \Big)\|\psi_1-\psi_2\|_{X_T}\\
        &\le  t^{\frac 14}\Big( t^{\frac 12}\|\n^4 u_1\|_{C^0}\|A\|_{C^1} +\|A\|_{C^0}+ t^{\frac 14}\|\n^3 u_1 \|_{C^0}+|\n^2 u_2 \|_{C^0}  \Big)\|\psi_1-\psi_2\|_{X_T}.
    \end{align*}
    Again, by definition of $\|\cdot\|_{X_T}$ and by \eqref{schauder Su},\eqref{higher order estimates S u} we conclude taking $ T,\ \delta$  small enough. 
\end{proof}

We are now able to prove a short-time existence result for the surface diffusion evolution. Thanks to the previous lemmas, we provide also higher order regularity estimates depending on the $C^{1,1}-$bound on the initial datum only. The proof follows closely the corresponding one in \cite{HZ,FusJulMor3D}.

\begin{theorem}\label{shortsurf}
    Let $\varepsilon>0$ and let $E \subset \T^N$ be a smooth open set.
    There exist $\delta=\delta(\varepsilon,E),$ $T=T(\varepsilon,E)>0$ with the following property: if $E_0$ is the normal deformation of $E$ induced by $u_0 \in C^{1,1}(\bd E)$, $\|u_0\|_{C^{1,1}(\bd E)}\le \delta$, and $|E_0|=|E|$, then the surface diffusion flow $E_t$ starting from $E_0$ exists in $[0,T)$, the sets $E_t$ are  normal deformations of $E$ induced by $u(\cdot,t) \in C^{\infty}(\bd E)$ for all $t \in (0,T)$, and 
    \begin{equation}\label{stima1surf}
        \sup_{t\in (0,T)} \|u\|_{C^{2}(\bd E)} \le \varepsilon.
    \end{equation}
    Moreover, for every $k \in \N \setminus \{0\}$, there exist constants $C_k=C_k(\varepsilon, E)>0$  such that
    \begin{equation}\label{stima2surf}
         \sup_{t\in [\frac{T}{2},T)} \|\n^{k+2} u \|_{C^0(\bd E)}\le C_{k} (\|u_0\|_{C^{1,1}(\bd E)}+1). 
    \end{equation}
\end{theorem}

\begin{proof}

In this proof we denote by $C>0$ a constant that depends on $N$ and $E$ and may change from line to line. Fix $\e>0$.\\
\textbf{Step 1:}  We show existence for \eqref{eq riscritta} via a fixed point argument. Let $T<1$, $\delta < 1$ to be chosen later, and let $u_1 \in C^{\infty}((0,T);C^\infty(\bd E))$ be the solution of 
\begin{equation*}
\begin{cases}
\bd _t u_1 =  -\Delta^2_E u_1 \quad &\text{on }\bd E\times [0,T), \vspace{5pt} \\
 u_1(\cdot, 0)=u_0  &\text{on }\bd E,
\end{cases}
\end{equation*}
where $u_0 \in C^{1,1}(\bd E)$ is such that $\|u_0\|_{C^{1,1}(\bd E)} \le \delta.$
The solution exists and it is given by \eqref{eqn4001}, that is $u_1=0+Su_0 \eqqcolon \psi_1 + Su_0$. Moreover \eqref{stima1surf} and \eqref{stima2surf} are satisfied by $u_1$ thanks to Theorem~\ref{thm001}, for $\delta$ small enough depending on $\e$. Let now $u_2$ be the solution of
\begin{equation*}
\begin{cases}
 \bd_t u_2=  -\Delta^2_E u_2 +f[u_{1}] \quad &\text{on }\bd E\times [0,T) ,\vspace{5pt}\\
u_2(\cdot, 0)=u_0\ &\text{on }\bd E,
\end{cases}
\end{equation*}
where $f[u]$ is defined as in \eqref{error term}. By \eqref{eqn4001} and \eqref{duhamel}, the unique solution  is given by $u_2=Vf[u_1]+Su_0=Vf[Su_0]+Su_0 \eqqcolon \psi_2 +Su_0.$ Moreover, by Theorem~\ref{thm701} and \eqref{est Jh} we have the estimate
\begin{equation*}
\begin{split}
    \|\psi_2\|_{X_T}\le C\|f[Su_0]\|_{Y_T}\le C\varepsilon(\|u_0\|_{C^{1,1}(\bd E)}+1) \le m,
\end{split}
\end{equation*}
for  $m$ sufficiently large. We are then led to define an iterative scheme. We set $u_1,u_2$ as above and for $n \ge 3$ we let $u_n$ be the solution to 
\begin{equation}
\begin{cases}\label{eq fixed point}
\bd_t u_n =  -\Delta^2_E u_n +f[u_{n-1}] \quad &\text{on }\bd E\times [0,T),\vspace{5pt}\\
u_n(\cdot, 0)=u_0\ &\text{on }\bd E,
\end{cases}
\end{equation}
and we split it as $u_n=S u_0+ V f[u_{n-1}]=:\psi_n +Su_0 .$ We will show that the sequence $\psi_n$ is converging in $X_T$. To do so, assume that  $\psi_j \in X_T$ for $j=1, \ldots, n-1$  with
\[\|\psi_j\|_{X_T}\le  m.\]
Then, by Theorem~\ref{thm701} and Lemma~\ref{lemma est Jh} we get $\psi_n \in X_T$ and
\begin{align}
    \|\psi_n\|_{X_T}& =\|V f[u_{n-1}]\|_{X_T}  \le C \|f[u_{n-1}]\|_{Y_T} = C \|f[\psi_{n-1}+Su_0]\|_{Y_T}\nonumber\\
    &\le C \sum_{j=2}^{n-1} \|f[\psi_{j}+S u_0]-f[\psi_{j-1}+S u_0]\|_{Y_T} +C\|f[Su_0]\|_{Y_T} \nonumber\\
    &\le C \Big (\sum_{j=1}^{n-1} \varepsilon^{j} \Big )(\|u_0\|_{C^{1,1}(\bd E)} + 1)  \nonumber\\
    &\le C \e\Big (1+\sum_{j=1}^{+\infty} \varepsilon^{j} \Big )(\|u_0\|_{C^{1,1}(\bd E)} + 1)  \nonumber\\
    &\le C \varepsilon (\|u_0\|_{C^{1,1}(\bd E)} + 1) \le m.\label{est X_T}
\end{align}
Moreover, Lemma \ref{lemma est Jh} implies that, for $\delta(\e,E), \ T(\e,E)$ small enough, it holds for all $n\geq 3$
\[  \| \psi_{n+1}-\psi_n \|_{X_T}\le \e \| \psi_{n}-\psi_{n-1} \|_{X_T},  \]
therefore $\psi_n$ is a Cauchy sequence and admits a limit point $\psi $ satisfying
\beq \label{bound fixed point}
\| \psi \|_{X_T}\le C\e (\|u_0\|_{C^{1,1}(\bd E)}+1).
\eeq
We thus showed the existence of a fixed point $u=\psi+ Su_0 $ for the problem \eqref{eq fixed point}. Finally, by \eqref{schauder Su} and \eqref{bound fixed point}  it holds
\beq\label{end step 1}
    \| u \|_{C^2(\bd E)}=\| \psi + S u_0\|_{C^2(\bd E)}\le \|\psi\|_{X_T}+ \|S u_0\|_{C^{2}(\bd E)}\le C\e (\| u_0\|_{C^{1,1}(\bd E)}+1).
\eeq

\medskip

\noindent\textbf{Step 2:} By \eqref{end step 1} we get straightforwardly that \eqref{stima2surf} holds for $k=0,1,2$.  In order to prove \eqref{stima2surf} for $k \ge 3$, we consider $x\in \bd E$ and we work under local coordinate, $B'_{r}\cong U \subset \bd E$ such that the metric $(g^{ij})_{i,j=1,\dots, N-1}$ of $\bd E$ satisfies $\frac 12 \delta_{ij}\le g_E^{ij}\le 2 \delta_{ij}$. Note in particular that  the operator $-\Delta^2_{E}$ is uniformly elliptic in $U$. In the following we identify $B'_r$ and $U\subset \bd E.$ We also set $g_t$ as the metric on $\bd E_t$ (see \cite[pag. 20]{Manlib} for details).  
Observe that $ u$ restricted to $B'_{r}\times [\frac{T}{2},T)$ is of class $C^{\infty}$ by the previous step. Recalling that $u=\psi+S u_0$, we have that the function $\psi$ satisfies 
\begin{equation}\label{eq psi}
\bd_t  \psi=-\Delta_{g_{t}}^2  \psi+ (\bd_t +\D^2_{g_{t}})(S u_0)+f'=: -\Delta_{g_{t}}^2  \psi +  \tilde f.
\end{equation}
Taking $\n_{g}$ in \eqref{eq psi} shows that the function $ \nabla_{g}\psi$ satisfies the equation
\begin{equation}
\begin{split}\label{eq psi deriv}
     \bd_t \nabla_{g} \psi&=-\Delta_{g_{t}}^2 \nabla_{g} \psi 
    -(\n_{g} g^{ij}_{t})g^{kl}_{t} (\psi)_{ijkl}-g^{ij}_{t} (\n_{g} g^{kl}_{t})(\psi)_{ijkl} +\n_{g}\tilde f
     \\
     &=:-\Delta_{g_{t}}^2 \nabla_{g} \psi  + F,
\end{split}
\end{equation}
where the error term $F$ contains the derivative of $\psi$ up to order four. 
To estimate $\|F\|_{C^{\beta/4}([\frac T2,T]; C^\beta(B'_r))}$ we first observe that, by \eqref{higher order estimates S u}, it follows
\begin{equation*}
\begin{split}
    \|\n_{g}& \left((\partial_t+\D^2_{g_{t}})(S {u_0})\right)\|_{C^{\beta/4}([\frac T2,T); C^\beta(B'_1))}\le  C \varepsilon( \|u_0\|_{C^{1,1}(\bd E)}+1).
    \end{split}
\end{equation*}
Secondly, we remark that the other terms of $F$ can be bounded analogously, recalling that they contain derivatives of $\psi$ up to order four and using \eqref{bound fixed point}, to show that
\begin{equation}\label{bound F lambda}
\begin{split}
\|F\|_{C^{\beta/4}([\frac T2,T); C^\beta(B'_r))} \le  C \varepsilon( \|u_0\|_{C^{1,1}(\bd E)}+1).
\end{split}
\end{equation}  
Note now that $\bd_t+\Delta^2_{g_{t}}$ is a uniformly parabolic operator, since the  coefficients of $\Delta^2_{g_t}$ are close to  the ones of $\Delta^2_E$  depending on $\| u(\cdot,t) \|_{C^{1,1}(\bd E)}$ as $g_{E_u}^{ij}-g_{E}^{ij} = B(x,u,\n u)$ and $B$ is a smooth function with $B(x,0,0)=0$, see again  \cite[pag. 20]{Manlib}. 
Since $\n_{g}\psi$ solves \eqref{eq psi deriv}, by the standard interior Schauder estimates and the bound \eqref{bound F lambda}, there exists $C>0$, which depends on $T$ and thus on $\varepsilon$ and $E$, such that
\begin{align*}
    \|\nabla_{g}\psi\|_{C^{1,\beta/4}([\frac T2, T); C^{4,\beta}(B'_{ r/2}))}&\le C\left( \|F\|_{C^{\beta/4}([\frac T4,T); C^\beta(B'_r))}   +\|\nabla_{g} \psi\|_{C^0(B'_r\times [\frac{T}{4},T))}\right)\\
    &\le C \varepsilon(\|u_0\|_{C^{1,1}(\bd E)} + 1),
\end{align*}
where we  noted  that $\|\psi\|_{C^1((B'_r\times [\frac{T}{4},T)))}\le \|\psi\|_{X_T}$ and employed again \eqref{bound fixed point}. Finally,
we conclude
\[\sup_{t \in [\frac T2,T)} \| \nabla^5  u\|_{C^0(\bd E)} \le C (\|u_0\|_{C^{1,1}(\bd E)} +1) .\]
By induction, one can prove \eqref{stima2surf} for every $k \in \N.$
\end{proof}

\section{Stability}

\subsection{Stability of the volume preserving mean curvature flow}
In this subsection, we study the evolution by mean curvature \eqref{defmcf} of normal deformations of a strictly stable set, as defined in Definition \ref{def strictly stable}. 
Suppose that $E$ is a strictly stable set and that $E_0=E_{u_0}$ is a smooth normal deformation of  $E$. By Theorem~\ref{teo short time MCF}, the volume preserving mean curvature flow starting from $E_0$ exists in a short time interval, and the evolving sets $E_t$ can be parametrized as normal deformations of the set $E$ induced by functions $u(\cdot,t)$ satisfying
\begin{equation*}
    \begin{cases}
    u_t(x,t) \nu_{E_t}(p)\cdot \nu_E(x)=-\left( \HHH_{E_t}(p)- \bar  \HHH_{E_t}\right) \qquad  x\in\bd E,\\
    u(\cdot,0)=u_0
    \end{cases}
\end{equation*}
where $p= x+u(x,t)\nu_E(x)$ and $\bar  \HHH_{E_t}=\fint_{\bd E_t}  \HHH_{E_t}  $. The scalar product above (see for instance \cite[eq.~(3.4)]{DeKu})  can be written as 
\[ \nu_{E_t}(p)\cdot \nu_E(x)=\left( 1+\sum_{j=1}^{N-1}\dfrac{(\bd_{\tau_j} u(x,t))^2}{(1+\kappa_j(x) u(x,t))^2} \right)^{ -1/2},\]
where $\kappa_j(x)$ and $\tau_j(x)$ are, respectively, the principal curvatures and the principal directions of $E$ at $x$. In particular, we remark that $\nu_{E_t}(p)\cdot \nu_E(x)=1+O(\|u(\cdot,t)\|_{H^1}).$  We can then prove the first part of the main result, that is Theorem~\ref{teo asymptotic}, concerning the long time behaviour of the volume preserving mean curvature flow. 

\begin{proof}[Proof of $(i)$ Theorem~\ref{teo asymptotic}]
Let $\varepsilon,\ \delta(\varepsilon) \in (0,1)$  to be chosen later. In the following, if not otherwise stated, the constants depend on $N,E$ and may change from line to line. Fix for instance $\beta= 1/2$ and suppose that $\delta$ is smaller than the constant given by Theorem~\ref{teo short time MCF}. We also use the short-hand notation $\pi_f\coloneqq(\pi_E|_{E_f})^{-1}$.\\
\textbf{Step 1.} We start by proving that $P(E_t)-P(E)\le C e^{-ct}$
as long as the flow exists.

Let $u_0\in C^{1,1}(\bd E)$ with $\|u_0\|_{C^{1,1}}\le \delta< 1$. By Theorem \ref{teo short time MCF} there exist a time $T>0 $, which depends on $E$ and the bound on $\|u_0\|_{C^{1,1}}< 1$, and a smooth flow $E_t$ starting from $E_0$ for $t\in[0,T)$. Moreover, $E_t=E_{u(\cdot,t)}$ and $u(\cdot,t)$ satisfies \eqref{stima1} and \eqref{stima2}. Without loss of generality we can assume $T<\infty$.
\\We notice that, considering $\e,\delta$ smaller, the value of $T$ does not change.

We recall the following well-known identities, holding along the smooth flow 
\begin{equation}\label{smooth comput}
    \dfrac \ud{\ud t}|E_t|=0,\quad \dfrac \ud{\ud t}P(E_t)= -\| \HHH_{E_t}-\bar \HHH_{E_t}  \|_{L^2(\bd E_t)}^2. 
\end{equation}
Let $\delta^*$ be the constant given by Theorem \ref{teo alex}, $p>N-1$ and $\eta=\eta(\delta^*, p)$ given by   Lemma~\ref{lemma 3.8_AFM}. By estimates \eqref{stima1}, \eqref{stima2} and by interpolation we have that  $\|u(\cdot,t)\|_{W^{2,p}(\bd E)}\le \eta$ for every $t\in [T/2,T)$, up to taking $\varepsilon$ smaller and therefore $\delta$ smaller. Thus for any $t\in [T/2,T)$ we can  apply Lemma~\ref{lemma 3.8_AFM} to find $\sigma_t\in \T^N$ and a function  $\tilde{u}(\cdot,t) $ such that $E_t +\sigma_t = E_{\tilde u(\cdot,t)}$ and
\begin{equation*}
\begin{split}
 |&\sigma_t| \le C \|u(\cdot,t)\|_{W^{2,p}(\bd E)},\ \|\tilde u(\cdot,t)\|_{W^{2,p}(\bd E)}\le C \|u(\cdot,t)\|_{W^{2,p}(\bd E)},\\
 &\left\lvert \int_{\bd E_t} \tilde u(\cdot,t)\nu_{E_t}  \right\rvert\le \delta^*\| \tilde u(\cdot,t) \|_{L^2(\bd E)}.
 \end{split}
 \end{equation*}
Furthermore, Lemma~\ref{lemma1.3} (taking $\delta$ smaller if needed) implies that $\| \tilde u(\cdot,t) \|_{C^1(\bd E)}\le \delta^*$. We then apply   Theorem~\ref{teo alex}  to the set $E_t+\sigma_t$ to obtain
\begin{equation}
    \| \tilde u(\cdot,t) \|_{H^1(\bd E)}\le C \| \mathscr H_{E_t+\sigma_t}-\lambda \|_{L^2(\bd E)}
\end{equation} 
for any $\lambda\in\R,$ where we recall $\mathscr H_{E_t+\sigma_t}(x)=\HHH_{E_t}(x+\tilde u(x)\nu_E(x)).$ From the previous equation, first by the change of variable  $y=x+\tilde u(x,t)\nu_E(x)$ (estimating the Jacobian with the bounds on $\tilde u$ and Lemma \ref{lemma1.3}), and then by translation invariance, we arrive at 
\begin{equation}\label{utildelesdifcurv}
    \| \tilde u(\cdot,t) \|_{H^1(\bd E)}\le C \| \HHH_{E_t+\sigma_t}-\lambda \|_{L^2(\bd E_t+\sigma_t)} =  C \| \HHH_{E_t}-\lambda \|_{L^2(\bd E_t)}.
\end{equation} 
We now claim that  
\begin{equation}\label{taylor perimeters}
     P(E_t+\sigma_t)-P(E)=P(E_{\tilde u(\cdot,t)})-P(E)\le C \|\tilde u(\cdot,t)\|^2_{H^1(\bd E)},
\end{equation}
which is a classical result but we provide a proof for the sake of completeness. 

Let us define, for every $x \in \bd E$, the function
$$
Q(x)\coloneqq \bigg( 1+ \sum_{j=1}^{N-1} \frac{(\partial_{\tau_j}\tilde{u}(x,t))^2}{(1+\kappa_j(x) \tilde{u}(x,t) )^2}   \bigg)^{\frac{1}{2}} $$
where $\tau_1(x),\ldots, \tau_{N-1}(x)$ and $\kappa_1(x),\ldots, \kappa_{N-1}(x)$ are, respectively, the principal directions and curvatures of $\bd E$ at $x$. Then by \cite[Lemma 3.1]{DeKu} we have
\begin{equation*}
\begin{split}
     P(E_t+ \sigma_t)& = P(E_{\tilde{u}(\cdot,t)}) = \int_{\partial E} Q(x)\prod_{i=1}^{N-1}\left( 1+ \kappa_i(x)\tilde{u}(t,\,x)\right)\,d \mathcal{H}^{N-1}(x) \\
    &=P(E)+ \int_{\partial E}\big(\HHH_{E}\tilde{u}(\cdot,t) + O(\tilde{u}(\cdot,t)^2)+ O(\vert D \tilde{u}(\cdot,t)\vert^2) \big)\, d \mathcal{H}^{N-1} \\
    & \leq P(E)+ \HHH_{E}\int_{\partial E} \tilde{u}(\cdot,t) \, d \mathcal{H}^{N-1}+ C\int_{\partial E} \big (\tilde{u}(\cdot,t)^2+ \vert D \tilde{u}(\cdot,t)\vert ^2 \big )\, d  \mathcal{H}^{N-1}\\
    &\le P(E)+ C \|\tilde{u}(\cdot,t) \|^2_{H^1(\partial E)},
    \end{split}
\end{equation*} 
where we have used that $ \HHH_E= \sum_{i=1}^{N-1} \kappa_i$ and the inequality $$ \left |\int_{\partial E} \tilde{u}(\cdot,t)\, d \mathcal{H}^{N-1}\right | \leq C \int_{\partial E} \tilde{u}(\cdot,t)^2 \, d \mathcal{ H}^{N-1},$$ which follows from the fact that $|E_t|=|E_0|$  (see \cite[Remark 3.2]{DeKu}). Hence, we prove the claim in \eqref{taylor perimeters}.

We now define the Lyapunov functional $\mathscr E(t)=P(E_t)-P(E)$, which is non increasing by \eqref{smooth comput}. Moreover, by translation invariance, from \eqref{utildelesdifcurv}, \eqref{taylor perimeters} and for any $\lambda\in\R$ we have
\begin{equation}\label{alex}
    \begin{split}
        P(E_t)-P(E)&=P(E_t+\sigma_t)-P(E)\le C \| \HHH_{E_t}-\lambda \|^2_{L^2(\bd E_t)}.
    \end{split}
\end{equation}
Since for any $t\in (0,T)$ equation  \eqref{alex} for the particular choice of $\lambda=\bar \HHH_{E_t} $ implies
\begin{align*}
    \mathscr E'(t) = -\|\HHH_{E_t}-\bar \HHH_{E_t}\|_{L^2(\bd E_t)}^2 \le -C \mathscr E (t),
\end{align*}
by Gronwall's inequality we conclude (recalling $\mathscr E(0)\ge \mathscr E( T/2)$)
\begin{equation}\label{decay}
    \mathscr E(t)\le \mathscr E(0) e^{-C(t-T/2)}, \qquad \forall t\in [T/2,T).
\end{equation}

\noindent\textbf{Step 2.} We now show that the flow exists for every $t\ge 0$ and it converges exponentially fast to $E$ up to translations.

Up to taking $\delta$ smaller, we can use the quantitative isoperimetric inequality in Theorem~\ref{coroll 1.2} to find the existence of  translations $\tau_t$ such that
\begin{equation*}
     C|E \triangle (E_t +\tau_t)|^2 \le P(E_t)-P(E) \le P(E_0)-P(E).
\end{equation*}
Furthermore, since all the evolving sets $\{ E_t\}_{t\in [T/2,T)}$ satisfy a uniform inner and outer ball condition by Remark \ref{rmk ball C11}, by classical convergence results (see e.g. \cite[Theorem 3.2]{Dal}) we have that $E_t+\tau_t$ is $C^{1}-$close to $E$.  In particular, there exist smooth (by the implicit map theorem) functions $v(\cdot,t): \bd E \to \R$ such that $E_t+\tau_t=E_{v(\cdot,t)}$ and 
$$|\tau_t| \le  \max_{x\in \bd E_t+\sigma_t} \dist_{\bd E_t}(x) \le  \|u(\cdot,t)\|_{C^0(\bd E)} +\|v(\cdot,t)\|_{C^{0}(\bd E)}\le 2 \e,$$
up to taking $\delta$ smaller. Therefore, recalling \eqref{decay}, we have
\begin{equation}\label{decay v}
    \|v(\cdot,t)\|_{L^1(\bd E)}^2\le C(P(E_0)-P(E)) e^{-C(t-T/2)}.
\end{equation}
By Lemma~\ref{lemma1.3}, we also have $\|v(\cdot,t)\|_{C^k(\bd E)}\le C( \|u(\cdot,t)\|_{C^k(\bd E)}+|\tau_t|)$ for every $k\ge 2$. 
For every $t\in [T/2,T)$, by combining the previous estimate with \eqref{stima2}, \eqref{decay v} and  interpolation inequalities, for any $l\in\N$ there exist $k(l)\in \N, \theta(l)\in (0,1)$ and $C=C(E,l)>0$ such that
\begin{equation}\label{decayC}
\|\nabla^l v(\cdot,t)\|_{C^{0}}\le C \|v(\cdot,t)\|_{L^1}^{\theta} \| v(\cdot,t) \|_{C^k}^{1-\theta}\le {C}{T^{-\frac{k}{4}(1-\theta)}}(P(E_0)-P(E))^{\frac{\theta}{2}} e^{-C (t-T/2)}.
\end{equation}

Choosing $\mathscr E(0)=P(E_0)-P(E)$ small (hence choosing $\delta$ small) we can then apply again Theorem~\ref{teo short time MCF} with the new initial set $E_{v(\cdot,T/2)}=E_{T/2} +\tau_{T/2}$ to get existence of the translated flow up to the time $3T/2.$  We remark that, by uniqueness, the flow above is well defined since it coincides in $[T/2,T)$ with the flow $E_t$ translated by $\tau_t$ and estimate \eqref{decay} now holds for all $t\in [T/2,3T/2).$   Since now the bound \eqref{decayC} is uniform along the flow, choosing at every step the times $t=n T/2$, we can iterate the procedure above to prove that the flow exists for all times $t\in[0,\infty)$. Moreover, for every $ t\in (0,\infty)$ there exists a translation $\tau_t$ such that $E_t+\tau_t=E_{v(\cdot,t)}$ with $v$ satisfying \eqref{decayC}. In particular, we have that $v\to 0$ exponentially in $C^k$ for any $k$,
as $t\to \infty$ and thus $E_t+\tau_t\to E$ in $C^k$ for every $k$. This also implies (reasoning as in \eqref{utildelesdifcurv}) that $\| \HHH_{E_t}-\bar \HHH_{E_t} \|_{L^2(\bd E)}\to 0$  exponentially fast.

\noindent\textbf{Step 3.} We conclude by showing the convergence of the whole flow to a translate of $E$.

Let us prove the convergence of the translations $\{\tau_t\}_{t\ge 0}$. 
By compactness we can find a sequence $t_n\to\infty$ such that  $\tau_{t_n} \to \tau$.  Defining 
\begin{equation}\label{def dissip}
\mathcal D(F,G)\coloneqq \int_{ F\triangle G}\dist_{\bd G}(x)\ud x,    
\end{equation}
following the computations of \cite[pag.~21]{AFJM} we see
\beq\label{decay D}
\begin{split}
     \left |\dfrac\ud{\ud t}\mathcal D(E_t, E-\tau)\right |&=\left |\dfrac\ud{\ud t} \int_{E_t\triangle (E-\tau)}\dist_{\bd E\tau_t}(x) \ud x\right |\\
     &=\left |\int_{E_t} \div(\sd_{E-\tau}(x) V_t(x) \nu_{E_t}(x)) \ud x\right |\\
     &=\left |-\int_{\bd E_t} \sd_{ E-\tau}(x) (\HHH_{E_t}(x)-\bar \HHH_{E_t}(x)) \ud \mathcal{H}^{N-1}(x)\right |\\
     &\le  P(E_0)\|\HHH_{E_t}-\bar \HHH_{E_t}\|_{L^2(\bd E)}  \left(\sup_{x \in \bd E_t}\dist_{\bd E-\tau}(x) \right) \\
     &\le C e^{-Ct} \left(\sup_{x \in \T^N}\dist_{\bd E-\tau}(x) \right)\le C  e^{-Ct},
\end{split}
\eeq
where we recall that $V_t$ is the velocity of the flow in the normal direction (see \eqref{defmcf}). Clearly, condition \eqref{decay D} implies that $\mathcal D(E_t, E-\tau)$ admits a limit as $t\to +\infty$. By the previous step and since $\tau_{t_n}\to \tau$, we deduce that 
\[ \mathcal D(E_t, E-\tau)\to 0 \quad \text{as }t\to+\infty. \]
Assume now that $\sigma\in\T^N$ is the limit of $\tau_{s_n}$ along a subsequence $s_n\to\infty$ as $n\to +\infty$. By the previous step, $E_{s_n}\to E-\sigma$, therefore 
\[  0=\lim_{n \rightarrow + \infty} \mathcal D(E_{s_n}, E-\tau)=\mathcal D(E-\sigma, E-\tau), \]
which implies $\sigma=\tau$ by definition \eqref{def dissip}. This concludes the proof as the exponential convergence follows from Step 2.
\end{proof}

\subsection{Stability of the surface diffusion flow}
We now focus on surface diffusion flow, which we defined in \eqref{defsdf}. 
As in the previous subsection, we consider $E$ a strictly stable set and $E_0=E_{u_0}$ a smooth normal deformation of  $E$. 
By Theorem~\ref{shortsurf}, the surface diffusion flow starting from $E_0$ exists smooth in an interval $[0,T)$, moreover the evolving sets $E_t$ can be written as normal deformations of  $E$ induced by functions $u(\cdot,t)$ satisfying
\begin{equation*}
    \begin{cases}
    u_t(x,t) \nu_{E_t}(p)\cdot \nu_E(x)=\Delta_{ E_t}\HHH_{E_t}(p)
    \qquad \forall x\in\bd E,\\
    u(x,0)=u_0(x)
    \end{cases}
\end{equation*}
where $p= x+u(x,t)\nu_E(x)$.

Now, we aim to show the stability result $(ii)$ of Theorem~\ref{teo asymptotic} for the surface diffusion flow. Due to the similarity of the arguments needed with those employed to prove item $(i)$ of Theorem~\ref{teo asymptotic}, we will only highlight the main differences between the two.

\begin{proof}[Proof of $(ii)$ Theorem~\ref{teo asymptotic}]
Firstly, Theorem~\ref{shortsurf} ensures the existence of a smooth flow $E_t$ for $t\in(0,T)$ of normal deformations of $E$ induced by functions $u(\cdot,t)\in C^{\infty}(\bd E)$ and satisfying \eqref{stima1surf} and \eqref{stima2surf}.
We recall the following identities, holding along the flow $E_t$ as long as it exists smooth,
\begin{equation}\label{smooth comput2}
    \dfrac \ud{\ud t}|E_t|=0,\quad \dfrac \ud{\ud t}P(E_t)=\int_{\partial E}\HHH_{E_t}(x)\Delta_{E_t}  \HHH_{E_t}(x)\, dx= -\|\nabla \HHH_{E_t}\|_{L^2(\bd E_t)}^2\le 0. 
\end{equation}
Denoting by $C_{E_t}$ the constant in the Poincaré inequality of Lemma~\ref{Poincaré inequality},  we get $$\| \HHH_{E_t}-\bar \HHH_{E_t} \|_{L^2(\bd E_t)}  \le  C_{E_t}\|\nabla \HHH_{E_t}\|_{L^2(\bd E_t)}.$$ 
Combining the previous inequality with \eqref{smooth comput2}, we obtain
\[ \dfrac \ud{\ud t}P(E_t)\le - C_{E_t}\| \HHH_{E_t}-\bar \HHH_{E_t} \|^2_{L^2(\bd E_t)}. \]
Since $\|u(\cdot,t)\|_{C^{1,1}(\bd E)}\le c$ for every $t\in (0,T),$ the Poincaré constants $C_{E_t}$ are uniformly bounded in the same time interval and the bound depends on $E,\|u\|_{C^{1,1}(\bd E)}$ (see e.g. the results in \cite{DDM}). Thus, we obtain the estimate $\frac{\ud}{\ud t}P(E_t)\le -C\|\HHH_{E_t}-\bar \HHH_{E_t}\|^2_{L^2(\bd E_t)}$ uniformly in $(0,T)$. 
We then conclude by following the same arguments of part $(i)$.
\end{proof}

\appendix
\section{Sketch of a general proof of the Lemma \ref{lemma est Jh} } 
In this appendix we complete the proof of Lemma \ref{lemma est Jh} in the general case, i.e. considering the full nonlinear error  term given by \eqref{formulaJestesa}.
\begin{proof}
As in Lemma \ref{lemma est Jh}, let $T<1$ to be chosen later.  We prove only equation \eqref{est Jh}, since the proof of the equations \eqref{contrazione}, \eqref{in space} is completely analogous.
\\We set
 \begin{equation*}
      f[u](x,t)\coloneqq \la A(x,u(x,t),\n u(x,t)), \n^4 u(x,t)  \ra + J(x,u (x,t),\nabla u (x,t),\nabla^2 u (x,t), \nabla^3 u (x,t)).
 \end{equation*}

Since the estimates for the first term of $f[u]$ have been presented in the proof of Lemma \ref{lemma est Jh}, we focus on bounding the terms of $J(x,u,\nabla u,\nabla^2 u, \nabla^3 u)$ with respect to the norm $\| \cdot \|_{C^0}$.
Considering the term $\langle B_1, \nabla^3 u \otimes \nabla^2 u \rangle $, we have
\begin{equation*}
 \langle B_1, \nabla^3 u \otimes \nabla^2 u \rangle \leq \| B_1 \|_{C^0} \| \nabla^3 u \otimes \nabla^2 u \|_{C^0} \leq C  \| \nabla^3 u \|_{C^0} \| \nabla^2 u \|_{C^0},
\end{equation*} 
as long as $\|u\|_{C^1}$ is small. Hence, with the same arguments presented for the functional $\la B, \n^3 u\otimes \n^2 u\ra$ we obtain
\begin{equation*}
\sup_{t \in (0,T)} t^{\frac{1}{2}}\langle B_1, \nabla^3 u \otimes \nabla^2 u \rangle  \leq \varepsilon \| u_0 \|_{C^{1,1}} 
\end{equation*} 
by choosing $T=T(\e)$ small enough. We analogously treat the other terms, so we have
\begin{equation*}
\langle B_2, \nabla^3 u\rangle \leq \| B_2 \|_{C^0} \|\nabla^3 u \|_{C^0} ,\, \langle B_3, \nabla^2 u \otimes \nabla^2 u \otimes \nabla^2 u \rangle \leq  \| B_3 \|_{C^0} \|\nabla^2 u \|_{C^0}^3  
\end{equation*}
and
\begin{equation*}
\langle B_4, \nabla^2 u \otimes \nabla^2 u \rangle \leq \| B_4 \|_{C^0} \| \nabla^2 u \|_{C^0}^2 ,\, \langle B_5, \nabla^2 u \rangle \leq \| B_5 \|_{C^0} \| \nabla^2 u \|_{C^0} \, .
\end{equation*}
Following the arguments of Lemma \ref{lemma est Jh}, we obtain
\begin{equation*}
\sup_{t \in (0,T)} t^{\frac{1}{2}} \langle B_2, \nabla^3 u\rangle  \leq \varepsilon \| u_0 \|_{C^{1,1}},
\end{equation*}
\begin{equation*} 
\sup_{t \in (0,T)} t^{\frac{1}{2}} \langle B_3, \nabla^2 u \otimes \nabla^2 u \otimes \nabla^2 u \rangle  \leq \varepsilon \| u_0 \|_{C^{1,1}},
\end{equation*}
\begin{equation*}
\sup_{t \in (0,T)} t^{\frac{1}{2}} \langle B_4, \nabla^2 u \otimes \nabla^2 u \rangle   \leq \varepsilon \| u_0 \|_{C^{1,1}}
\end{equation*}
and
\begin{equation*}
 \sup_{t \in (0,T)} t^{\frac{1}{2}} \langle B_5, \nabla^2 u \rangle \leq \varepsilon \| u_0 \|_{C^{1,1}}.
\end{equation*}

In the end we have that $\|b_6\|_{C^0}\le C=C(E)$. 
\\Therefore, taking $T$ small we obtain 
\[  \sup_t  t^{\frac 12}\|b_6\|_{C^0}\le\e.\]
We now focus on the H\"older seminorm in space. We present the calculations only for $\langle B_1, \nabla^3 u \otimes \nabla^2 u \rangle$, being the other analogous. A straightforward computation shows (using the triangular inequality) that
\begin{equation*} \begin{split}
 & | \langle B_1(x+h, u(x+h), \n u(x+h)), \nabla^3 u (x+h)\otimes \nabla^2 u(x+h) \rangle \\ &- \langle B_1 (x, u(x), \n u(x)), \nabla^3 u (x)\otimes \nabla^2 u(x) \rangle | \\
 & \leq (|h | \| \partial_1 B_1\|_{C^0}+ | u(x+h)-u(x)| \| \partial_2 B_1 \|  \\ & + | \n u(x+h)- \n u(x) | \| \partial_3 B_1 \|_{C^0}  ) \| \n^3 u\|_{C^0}\| \n^2 u \|_{C^0}   \\
 & + \| B_1 \|_{C^0} ( |\n^3(x+h)-\n^3 u(x) | \| \n^2 u \|_{C^0}+ \|\n^2 u\|_{C^0} | \n^2 u(x+h)- \n^2 u(x) |).
\end{split}
\end{equation*}
Therefore, as in the case $J(x,u,\n u, \n^2 u, \n^3 u)= \la B, \n^3 u \times \n^2 u\ra$, using formula \eqref{schauder Su} and \eqref{higher order estimates S u} we obtain the thesis.

Finally, we show how to bound  the H\"older seminorm in time appearing in $\|f[u]\|_{Y_T}$. We fix $t\in (0,T),\tau\in (0,T-t)$ and, for simplicity, we omit the dependence on $x$. 
For the first term, we have
\begin{equation*}
\begin{split}
   &| \langle B_1, \nabla^3 u (t+\tau) \otimes \nabla^2 u (t+\tau)\rangle - \langle B_1, \nabla^3 u (t)\otimes \nabla^2 u (t) \rangle | \\ & \leq \| B_1 \|_{C^0} \left[ | \n^3 u(t+\tau)- \n^3 u(t) | \|\nabla^2 u \|_{C^0}+ | \n^2 u(t+\tau)- \n^2 u(t)| \|\n^3 u \|_{C^0}\right] \\
   & +\big(\| \partial_2 B_1 \|_{C^0} | u(t)-u(t+\tau) |+ \| \partial_3 B_1 \| | \nabla u(t)- \nabla u(t+\tau) | \big) \| \n^3 u\|_{C^0} \| \n^2 u\|_{C^0}.
   \end{split}
\end{equation*} 
Then, for the second, third and fourth terms we get, respectively,
\begin{equation*}
\begin{split}
  | \langle B_2,& \nabla^3 u(t+\tau)\rangle- \langle B_2, \nabla^3 u(t)\rangle|  \leq \| B_2 \|_{C^0} | \n^3 u(t+\tau)- \n^3 u(t) |  \\
  &+\big(\| \partial_2 B_2 \|_{C^0} | u(t)-u(t+\tau) |+ \| \partial_3 B_2 \| | \nabla u(t)- \nabla u(t+\tau) | \big) \| \n^3 u\|_{C^0},
  \end{split}
\end{equation*} 
\begin{equation*}
 \begin{split}
       & | \langle B_3, \nabla^2 u (t+\tau) \otimes \nabla^2 u (t+\tau) \otimes \nabla^2 u (t+\tau) \rangle -\langle B_3, \nabla^2 u (t)\otimes \nabla^2 u (t) \otimes \nabla^2 u (t) \rangle  | \\
       & \leq 3 \| B_3\|_{C^0} | \n^2 u(t+\tau)- \n^2 u(t) | \| \n^2 u \|_{C^0}^2\\
  &\quad+\big(\| \partial_2 B_3 \|_{C^0} | u(t)-u(t+\tau) |+ \| \partial_3 B_3 \| | \nabla u(t)- \nabla u(t+\tau) | \big) \| \n^2 u\|_{C^0}^3,
     \end{split}
\end{equation*}
and
\begin{equation*}    \begin{split}
     & |\langle B_4, \nabla^2 u(t+\tau) \otimes \nabla^2 u(t+\tau) \rangle - \langle B_4, \nabla^2 u(t) \otimes \nabla^2 u(t) \rangle| \\
     &  \leq 2 \|B_4 \|_{C^0} \|\n^2 u \|_{C^0} |\n^2 u(t+\tau)-\n^2 u(t) |\\
  &\quad+\big(\| \partial_2 B_4 \|_{C^0} | u(t)-u(t+\tau) |+ \| \partial_3 B_4 \| | \nabla u(t)- \nabla u(t+\tau) | \big) \| \n^2 u\|_{C^0}^2.
 \end{split}
\end{equation*} 
Finally, for the last two terms we have
\begin{equation*}
\begin{split}
 | \langle& B_5, \nabla^2 u (t+\tau)\rangle - \langle B_5, \nabla^2 u (t)\rangle | \leq \| B_5 \|_{C^0} | \n^2 u(t+\tau)- \nabla^2 u(t)|\\
  &+\big(\| \partial_2 B_5 \|_{C^0} | u(t)-u(t+\tau) |+ \| \partial_3 B_5 \| | \nabla u(t)- \nabla u(t+\tau) | \big) \| \n^2 u\|_{C^0}
  \end{split}
\end{equation*}
and
\begin{equation*}
 \begin{split}
       &  | b_6(\cdot, u(t+\tau), \n u(t+\tau))- b_6(\cdot, u(t),\n u(t)) | \\ &\leq  (\| \partial_2 b_6 \|_{C^0}+ \| \partial_3 b_6 \|_{C^0})( | u(t+\tau)-u(t) |+ | \n u(t+\tau)-\n u(t) |).
     \end{split}
\end{equation*}
Therefore, we can conclude with the same arguments  used for $\la B, \n^3 u\otimes \n^2 u\ra$.
\end{proof}

\section*{Conflict of interest}
The authors declare that they have no conflict of interest.

\section*{Acknowledgements}

The authors warmly thank Massimiliano Morini for many helpful discussions and advice.
Antonia Diana and Anna Kubin are members of the Gruppo Nazionale per l’Analisi Matematica, la Probabilità e le loro Applicazioni (GNAMPA) of the Istituto Nazionale di Alta Matematica (INdAM).  
Daniele De Gennaro has received funding from the European Union’s Horizon 2020 research and innovation programme under the Marie Skłodowska-Curie grant agreement No 94532.
Andrea Kubin is supported by the DFG Collaborative Research Center TRR 109 “Discretization in Geometry and Dynamics”. The research of Anna Kubin was funded in whole or in part by the INdAM--GNAMPA 2023 Project \textit{Problemi variazionali per funzionali e operatori non-locali} (grant agreement No. CUP\_E53\-C22\-001\-930\-001), and by the Austrian Science Fund (FWF) [10.55776/F65] and [10.55776/P35359].

\printbibliography 

@article {GarckGob1,
    AUTHOR = {Garcke, H. and G\"{o}\ss wein, M.},
     TITLE = {On the surface diffusion flow with triple junctions in higher
              space dimensions},
   JOURNAL = {Geom. Flows},
  FJOURNAL = {Geometric Flows},
    VOLUME = {5},
      YEAR = {2020},
    NUMBER = {1},
     PAGES = {1--39},
       DOI = {10.1515/geofl-2020-0001},
}

@article{GarckGob2,
	author = {H. Garcke and M. G{\"o}{\ss}wein},
	doi = {https://doi.org/10.1016/j.jde.2021.08.033},
	issn = {0022-0396},
	journal = {Journal of Differential Equations},
	pages = {617-661},
	title = {Non-linear stability of double bubbles under surface diffusion},
	volume = {302},
	year = {2021},
}

@book{Manlib,
	author = {C. Mantegazza},
	publisher = {Birkh\"auser/Springer Basel AG, Basel},
	series = {Progress in Mathematics},
	title = {Lecture notes on mean curvature flow},
	volume = {290},
	year = {2011}}

@article {DDMsurvey,
    AUTHOR = {Della Corte, S. and Diana, A. and Mantegazza, C.},
     TITLE = {Global existence and stability for the modified
              {M}ullins-{S}ekerka and surface diffusion flow},
   JOURNAL = {Math. Eng.},
  FJOURNAL = {Mathematics in Engineering},
    VOLUME = {4},
      YEAR = {2022},
    NUMBER = {6},
     PAGES = {Paper No. 054, 104},
   MRCLASS = {53E10},
  MRNUMBER = {4354990},
       DOI = {10.3934/mine.2022054},
}

@unpublished{DFM,
	author = {A. Diana and N. Fusco and C. Mantegazza},
	note = {arXiv:2304.04011v1},
	title = {{S}tability for the {S}urface {D}iffusion {F}low},
	year={2023}}

@unpublished{DDM,
	author = {S. {Della Corte} and A. Diana and C. Mantegazza},
	note = {{arXiv:2107.12234v4}},
	title = {Uniform {S}obolev, interpolation and {Calder\'on--Zygmund} inequalities for graph hypersurfaces},
	year={2023}
 %URL={https://doi.org/10.48550/arXiv.2304.04013}
 }

@book{gilbarg1977elliptic,
  title={Elliptic partial differential equations of second order},
  author={Gilbarg, D. and Trudinger, N.},
  volume={224},
  number={2},
  year={1977},
  publisher={Springer}
}

@article {JulMorPonSpa,
    AUTHOR = {Julin, V. and Morini, M. and Ponsiglione, M. and Spadaro, E.},
     TITLE = {The asymptotics of the area-preserving mean curvature and the Mullins–Sekerka flow in two dimensions},
   JOURNAL = {Math. Ann.},
  FJOURNAL = {Mathematische Annalen},
  YEAR = {2022},
       DOI = {10.1007/s00208-022-02497-3},
}

@article {Dal,
    AUTHOR = {Dalphin, J.},
     TITLE = {Uniform ball property and existence of optimal shapes for a
              wide class of geometric functionals},
   JOURNAL = {Interfaces Free Bound.},
  FJOURNAL = {Interfaces and Free Boundaries. Mathematical Analysis,
              Computation and Applications},
    VOLUME = {20},
      YEAR = {2018},
    NUMBER = {2},
     PAGES = {211--260},
      %ISSN = {1463-9963},
       DOI = {10.4171/IFB/401},
}

@incollection {Amalib,
    AUTHOR = {Amann, H.},
     TITLE = {Nonhomogeneous linear and quasilinear elliptic and parabolic
              boundary value problems},
 BOOKTITLE = {Function spaces, differential operators and nonlinear analysis
              ({F}riedrichroda, 1992)},
    SERIES = {Teubner-Texte Math.},
    VOLUME = {133},
     PAGES = {9--126},
 PUBLISHER = {Teubner, Stuttgart},
      YEAR = {1993},
   MRCLASS = {35K57 (34G20 35J60 46E35 47N20)},
  MRNUMBER = {1242579},
       DOI = {10.1007/978-3-663-11336-2\_1},
       URL = {https://doi.org/10.1007/978-3-663-11336-2_1},
}

@article {DegKuKu,
    AUTHOR = {De Gennaro, D. and Kubin, A. and Kubin, A.},
     TITLE = {Asymptotic of the discrete volume-preserving fractional mean
              curvature flow via a nonlocal quantitative {A}lexandrov
              theorem},
   JOURNAL = {Nonlinear Anal.},
  FJOURNAL = {Nonlinear Analysis. Theory, Methods \& Applications. An
              International Multidisciplinary Journal},
    VOLUME = {228},
      YEAR = {2023},
     PAGES = {Paper No. 113200, 23},
      %ISSN = {0362-546X},
   MRCLASS = {35R11 (49M25 49Q20 53E10)},
  MRNUMBER = {4523489},
       DOI = {10.1016/j.na.2022.113200},
}

@incollection {BBM,
    AUTHOR = {Bourgain, J. and Brezis, H. and Mironescu, P.},
     TITLE = {Another look at {S}obolev spaces},
 BOOKTITLE = {Optimal control and partial differential equations},
     PAGES = {439--455},
 PUBLISHER = {IOS, Amsterdam},
      YEAR = {2001},
   MRCLASS = {46E30},
  MRNUMBER = {3586796},
}

@article {ES,
    AUTHOR = {Escher, J. and Simonett, G.},
     TITLE = {The volume preserving mean curvature flow near spheres},
   JOURNAL = {Proc. Amer. Math. Soc.},
  FJOURNAL = {Proceedings of the American Mathematical Society},
    VOLUME = {126},
      YEAR = {1998},
    NUMBER = {9},
     PAGES = {2789--2796},
   MRCLASS = {53C21 (35B40 35K55 53A07 53C40 58G11)},
  MRNUMBER = {1485470},
MRREVIEWER = {Barbara Priwitzer},
       DOI = {10.1090/S0002-9939-98-04727-3},
}

@article {EMS,
    AUTHOR = {Escher, J. and Mayer, U. F. and Simonett, G.},
     TITLE = {The surface diffusion flow for immersed hypersurfaces},
   JOURNAL = {SIAM J. Math. Anal.},
  FJOURNAL = {SIAM Journal on Mathematical Analysis},
    VOLUME = {29},
      YEAR = {1998},
    NUMBER = {6},
     PAGES = {1419--1433},
   MRCLASS = {58E15 (35K99 35R35 65C20)},
  MRNUMBER = {1638074},
MRREVIEWER = {Vladimir Grebenev},
       DOI = {10.1137/S0036141097320675},
}

@article {Li,
    AUTHOR = {Li, H.},
     TITLE = {The volume-preserving mean curvature flow in {E}uclidean
              space},
   JOURNAL = {Pacific J. Math.},
  FJOURNAL = {Pacific Journal of Mathematics},
    VOLUME = {243},
      YEAR = {2009},
    NUMBER = {2},
     PAGES = {331--355},
      %ISSN = {0030-8730},
   MRCLASS = {53C44},
  MRNUMBER = {2552262},
MRREVIEWER = {James Alexander McCoy},
       DOI = {10.2140/pjm.2009.243.331},
}

@article {Whe,
    AUTHOR = {Wheeler, G.},
     TITLE = {Surface diffusion flow near spheres},
   JOURNAL = {Calc. Var. Partial Differential Equations},
  FJOURNAL = {Calculus of Variations and Partial Differential Equations},
    VOLUME = {44},
      YEAR = {2012},
    NUMBER = {1-2},
     PAGES = {131--151},
      %ISSN = {0944-2669},
   MRCLASS = {53C44 (58J35)},
  MRNUMBER = {2898774},
MRREVIEWER = {Weiyong He},
       DOI = {10.1007/s00526-011-0429-4},
}

@book{PruSim,
    author ={Prüss, J. and Simonett, G.},
    title = {Moving interfaces and quasilinear parabolic evolution equations},
    notes={Monographs in Mathematics, 105},
    publisher ={Birkhäuser/Springer} ,
    pages={xix+609},
    year = {2016}
}

@article {LeCSim20,
    AUTHOR = {LeCrone, J. and Simonett, G.},
     TITLE = {On quasilinear parabolic equations and continuous maximal regularity},
   JOURNAL = {Evol. Equ. Control Theory},
  FJOURNAL = {},
    VOLUME = {9},
      YEAR = {2020},
    NUMBER = {1},
     PAGES = {61-86},
      %ISSN = {0036-1410},
   MRCLASS = {35K93 (35B30 35B32 35B35 53C44)},
  MRNUMBER = {3103249},
MRREVIEWER = {Glen E. Wheeler},
       DOI = {10.1137/120883505},
}

@article {LeCSim13,
    AUTHOR = {LeCrone, J. and Simonett, G.},
     TITLE = {On well-posedness, stability, and bifurcation for the
              axisymmetric surface diffusion flow},
   JOURNAL = {SIAM J. Math. Anal.},
  FJOURNAL = {SIAM Journal on Mathematical Analysis},
    VOLUME = {45},
      YEAR = {2013},
    NUMBER = {5},
     PAGES = {2834--2869},
      %ISSN = {0036-1410},
   MRCLASS = {35K93 (35B30 35B32 35B35 53C44)},
  MRNUMBER = {3103249},
MRREVIEWER = {Glen E. Wheeler},
       DOI = {10.1137/120883505},
}

@article {AbeAraGar,
    AUTHOR = {Abels, H. and Arab, N. and Garcke, H.},
     TITLE = {Standard planar double bubbles are dynamically stable under
              surface diffusion flow},
   JOURNAL = {Comm. Anal. Geom.},
  FJOURNAL = {Communications in Analysis and Geometry},
    VOLUME = {29},
      YEAR = {2021},
    NUMBER = {5},
     PAGES = {1007--1060},
   MRCLASS = {53E10 (58J65)},
  MRNUMBER = {4349138},
MRREVIEWER = {James Alexander McCoy},
       DOI = {10.4310/CAG.2021.v29.n5.a1},
}

@article {GarItoKoh,
    AUTHOR = {Garcke, H. and Ito, K. and Kohsaka, Y.},
     TITLE = {Nonlinear stability of stationary solutions for surface
              diffusion with boundary conditions},
   JOURNAL = {SIAM J. Math. Anal.},
  FJOURNAL = {SIAM Journal on Mathematical Analysis},
    VOLUME = {40},
      YEAR = {2008},
    NUMBER = {2},
     PAGES = {491--515},
      %ISSN = {0036-1410},
   MRCLASS = {35K55 (35B35 35B45 35R35 53C44 74N20)},
  MRNUMBER = {2438774},
MRREVIEWER = {Jana Kopfova},
       DOI = {10.1137/070694752},
}

@incollection {Gag,
    AUTHOR = {Gage, M.},
     TITLE = {On an area-preserving evolution equation for plane curves},
 BOOKTITLE = {Nonlinear problems in geometry ({M}obile, {A}la., 1985)},
    SERIES = {Contemp. Math.},
    VOLUME = {51},
     PAGES = {51--62},
 PUBLISHER = {Amer. Math. Soc., Providence, RI},
      YEAR = {1986},
   MRCLASS = {53A04},
  MRNUMBER = {848933},
MRREVIEWER = {W. J. Firey},
       DOI = {10.1090/conm/051/848933},
}

@article {FusJulMor3D,
    AUTHOR = {Fusco, N. and Julin, V. and Morini, M.},
     TITLE = {The surface diffusion flow with elasticity in three dimensions},
   JOURNAL = {Arch. Ration. Mech. Anal.},
  FJOURNAL = {Archive for Rational Mechanics and Analysis},
    VOLUME = {237},
      YEAR = {2020},
    NUMBER = {3},
     PAGES = {1325--1382},
      %ISSN = {0003-9527},
   MRCLASS = {74G65 (49Q10)},
  MRNUMBER = {4110437},
MRREVIEWER = {Elvira Zappale},
       DOI = {10.1007/s00205-020-01532-4},
}

@article {FusJulMor2D,
    AUTHOR = {Fusco, N. and Julin, V. and Morini, M.},
     TITLE = {The surface diffusion flow with elasticity in the plane},
   JOURNAL = {Comm. Math. Phys.},
  FJOURNAL = {Communications in Mathematical Physics},
    VOLUME = {362},
      YEAR = {2018},
    NUMBER = {2},
      %ISSN = {0010-3616},
   MRCLASS = {35Q74 (74N25)},
  MRNUMBER = {3843423},
MRREVIEWER = {Giuliano Lazzaroni},
       DOI = {10.1007/s00220-018-3200-2},
}

@article {MaySim,
    AUTHOR = {Mayer, U. F. and Simonett, G.},
     TITLE = {Self-intersections for the surface diffusion and the
              volume-preserving mean curvature flow},
   JOURNAL = {Differential Integral Equations},
  FJOURNAL = {Differential and Integral Equations. An International Journal
              for Theory \& Applications},
    VOLUME = {13},
      YEAR = {2000},
    NUMBER = {7-9},
     PAGES = {1189--1199},
      %ISSN = {0893-4983},
}

@article {Mul,
    AUTHOR = {Mullins, W. W.},
     TITLE = {Theory of Thermal Grooving},
   JOURNAL = {Journal of Applied Physics},
    VOLUME = {28},
      YEAR = {1957},
     PAGES = {333--339},
}

@article {AFM,
    AUTHOR = {Acerbi, E. and Fusco, N. and Morini, M.},
     TITLE = {Minimality via second variation for a nonlocal isoperimetric
              problem},
   JOURNAL = {Comm. Math. Phys.},
  FJOURNAL = {Communications in Mathematical Physics},
    VOLUME = {322},
      YEAR = {2013},
    NUMBER = {2},
     PAGES = {515--557},
   MRCLASS = {49Q05 (58E30)},
  MRNUMBER = {3077924},
MRREVIEWER = {Martin Fuchs},
       DOI = {10.1007/s00220-013-1733-y},
}

@article {MorPonSpa,
    AUTHOR = {Morini, M. and Ponsiglione, M. and Spadaro, E.},
     TITLE = {Long time behavior of discrete volume preserving mean
              curvature flows},
   JOURNAL = {J. Reine Angew. Math.},
  FJOURNAL = {Journal f\"{u}r die Reine und Angewandte Mathematik. [Crelle's
              Journal]},
    VOLUME = {784},
      YEAR = {2022},
     PAGES = {27--51},
      ISSN = {0075-4102},
   MRCLASS = {53E10},
  MRNUMBER = {4388339},
MRREVIEWER = {Yue He},
       DOI = {10.1515/crelle-2021-0076},
       URL = {https://doi.org/10.1515/crelle-2021-0076},
}

@article {BelCasChaNov,
    AUTHOR = {Bellettini, G. and Caselles, V. and Chambolle,
              A. and Novaga, M.},
     TITLE = {Crystalline mean curvature flow of convex sets},
   JOURNAL = {Arch. Ration. Mech. Anal.},
  FJOURNAL = {Archive for Rational Mechanics and Analysis},
    VOLUME = {179},
      YEAR = {2006},
    NUMBER = {1},
     PAGES = {109--152},
   MRCLASS = {53C44 (58E12)},
  MRNUMBER = {2208291},
MRREVIEWER = {Matthias Wilhelm Winter},
       DOI = {10.1007/s00205-005-0387-0},
}

@article {DeKu,
    AUTHOR = {De Gennaro, D. and Kubin, A.},
     TITLE = {Long time behaviour of the discrete volume preserving mean
              curvature flow in the flat torus},
   JOURNAL = {Calc. Var. Partial Differential Equations},
  FJOURNAL = {Calculus of Variations and Partial Differential Equations},
    VOLUME = {62},
      YEAR = {2023},
    NUMBER = {3},
     PAGES = {Paper No. 103, 39},
   MRCLASS = {49M25 (49Q20 53C24 53E10)},
  MRNUMBER = {4553954},
       DOI = {10.1007/s00526-023-02439-0},
}

@article {Nii,
    AUTHOR = {Niinikoski, J.},
     TITLE = {Volume preserving mean curvature flows near strictly stable
              sets in flat torus},
   JOURNAL = {J. Differential Equations},
  FJOURNAL = {Journal of Differential Equations},
    VOLUME = {276},
      YEAR = {2021},
     PAGES = {149--186},
      %ISSN = {0022-0396},
   MRCLASS = {53E10 (35K93)},
  MRNUMBER = {4191351},
MRREVIEWER = {Kazuhiro Horihata},
       DOI = {10.1016/j.jde.2020.12.010},
}

@article {JN,
    AUTHOR = {Julin, V. and Niinikoski, J.},
     TITLE = {Consistency of the flat flow solution to the volume preserving
              mean curvature flow},
   JOURNAL = {Arch. Ration. Mech. Anal.},
  FJOURNAL = {Archive for Rational Mechanics and Analysis},
    VOLUME = {248},
      YEAR = {2024},
    NUMBER = {1},
     PAGES = {Paper No. 1, 58},
       DOI = {10.1007/s00205-023-01944-y},
}

@article{Mullins,
	author = {Mullins, W. W.},
	doi = {10.1063/1.1722511},
	fjournal = {Journal of Applied Physics},
	ISSN = {0944-2669},
	journal = {Journal of Applied Physics},
	number = {},
	pages = {900},
	title = {Two‐Dimensional Motion of Idealized Grain Boundaries},
	volume = {27},
	year = {1956},
}

@article {LS,
    AUTHOR = {Luckhaus, S. and Sturzenhecker, T.},
     TITLE = {Implicit time discretization for the mean curvature flow
              equation},
   JOURNAL = {Calc. Var. Partial Differential Equations},
  FJOURNAL = {Calculus of Variations and Partial Differential Equations},
    VOLUME = {3},
      YEAR = {1995},
    NUMBER = {2},
     PAGES = {253--271},
      %ISSN = {0944-2669},
   MRCLASS = {65M06 (58E12)},
  MRNUMBER = {1386964},
MRREVIEWER = {R. N. Mukherjee},
       DOI = {10.1007/BF01205007},
       URL = {https://doi.org/10.1007/BF01205007},
}

@article {Hui,
    AUTHOR = {Huisken, G.},
     TITLE = {The volume preserving mean curvature flow},
   JOURNAL = {J. Reine Angew. Math.},
  FJOURNAL = {Journal f\"{u}r die Reine und Angewandte Mathematik. [Crelle's
              Journal]},
    VOLUME = {382},
      YEAR = {1987},
     PAGES = {35--48},
      %ISSN = {0075-4102},
   MRCLASS = {53A07 (35K99 53C40)},
  MRNUMBER = {921165},
MRREVIEWER = {Michael Gr\"{u}ter},
       DOI = {10.1515/crll.1987.382.35},
}

@article {ATW,
    AUTHOR = {Almgren, F. and Taylor, J. E. and Wang, L.},
     TITLE = {Curvature-driven flows: a variational approach},
   JOURNAL = {SIAM J. Control Optim.},
  FJOURNAL = {SIAM Journal on Control and Optimization},
    VOLUME = {31},
      YEAR = {1993},
    NUMBER = {2},
     PAGES = {387--438},
      %ISSN = {0363-0129},
   MRCLASS = {58E50 (49Q10 58E15 73B99)},
  MRNUMBER = {1205983},
       DOI = {10.1137/0331020},
       URL = {https://doi.org/10.1137/0331020},
}

@book {Mag,
    AUTHOR = {Maggi, F.},
     TITLE = {Sets of finite perimeter and geometric variational problems},
    SERIES = {Cambridge Studies in Advanced Mathematics},
    VOLUME = {135},
      NOTE = {An introduction to geometric measure theory},
 PUBLISHER = {Cambridge University Press, Cambridge},
      YEAR = {2012},
     PAGES = {xx+454},
       DOI = {10.1017/CBO9781139108133},
}

@article {KocLam,
    AUTHOR = {Koch, H. and Lamm, T.},
     TITLE = {Geometric flows with rough initial data},
   JOURNAL = {Asian J. Math.},
  FJOURNAL = {Asian Journal of Mathematics},
    VOLUME = {16},
      YEAR = {2012},
    NUMBER = {2},
     PAGES = {209--235},
      %ISSN = {1093-6106},
   MRCLASS = {53C44},
  MRNUMBER = {2916362},
MRREVIEWER = {Meng Zhu},
       DOI = {10.4310/AJM.2012.v16.n2.a3},
}

@book {Fri,
    AUTHOR = {Friedman, A.},
     TITLE = {Partial differential equations of parabolic type},
 PUBLISHER = {Prentice-Hall, Inc., Englewood Cliffs, N.J.},
      YEAR = {1964},
     PAGES = {xiv+347},
   MRCLASS = {35.00 (35.62)},
  MRNUMBER = {0181836},
MRREVIEWER = {B. Frank Jones, Jr.},
}

@article {HZ,
    AUTHOR = {He, W. and Zeng, Y.},
     TITLE = {The {C}alabi flow with rough initial data},
   JOURNAL = {Int. Math. Res. Not. IMRN},
  FJOURNAL = {International Mathematics Research Notices. IMRN},
      YEAR = {2021},
    NUMBER = {10},
      %ISSN = {1073-7928},
   MRCLASS = {53E30 (32Q15)},
  MRNUMBER = {4259154},
MRREVIEWER = {Masaya Kawamura},
       DOI = {10.1093/imrn/rnz050},
}

@article {AFJM,
    AUTHOR = {Acerbi, E. and Fusco, N. and Julin, V. and Morini, M.},
     TITLE = {Nonlinear stability results for the modified {M}ullins-{S}ekerka and the surface diffusion flow},
   JOURNAL = {J. Differential Geom.},
  FJOURNAL = {Journal of Differential Geometry},
    VOLUME = {113},
      YEAR = {2019},
    NUMBER = {1},
     PAGES = {1--53},
     DOI = {10.4310/jdg/1567216953},
}

@book {Aub,
    AUTHOR = {Aubin, T.},
     TITLE = {Some nonlinear problems in {R}iemannian geometry},
    SERIES = {Springer Monographs in Mathematics},
 PUBLISHER = {Springer-Verlag, Berlin},
      YEAR = {1998},
     PAGES = {xviii+395},
       DOI = {10.1007/978-3-662-13006-3},
}
\end{document}